\documentclass[preprint,11pt]{elsarticle}
\usepackage[paperwidth=7.5in, paperheight=11in, margin=.775in]{geometry}
\usepackage{amsfonts}
\usepackage{amsmath}
\usepackage{graphicx}
\usepackage{graphics}
\usepackage{epsfig, setspace}
\usepackage{mathrsfs}
\usepackage{amssymb}
\usepackage{amsfonts}
\usepackage[normalem]{ulem}
\usepackage{xcolor}
\usepackage{todonotes}

\newcommand{\removed}[1] {\ifmmode{\color{red}\cancel{ #1}}\else{\color{red}\sout{#1}}\fi}

\newcommand{\tnorm}[1]{\ensuremath{\left| \! \left| \! \left|} #1\ensuremath{\right| \! \right| \! \right|}}

\def\cT{\mathcal{K}}

\def\cM{\mathcal{M}}
\def\cT{\mathcal{T}}
\def\T{\mathcal{T}}
\def\cK{\mathcal{T}}
\def\cV{\mathcal{V}}

\def\cE{\mathcal{E}}

\def\cE{\mathcal{E}}
\def\E{\mathcal{E}}

\def\K{\mathcal{T}}
\newcommand{\dx}{\;\textit{dx}}
\newcommand{\ds}{\;\textit{ds}}

\newcommand{\norm}[1]{\left\Vert #1 \right\Vert}

\newtheorem{theorem}{Theorem}[section]
\newtheorem{lemma}{Lemma}[section]

\begin{document}

\title{An overlapping local projection stabilization for Galerkin approximations of Darcy and Stokes problems}
%\bibliographystyle{amsplain.bst}
%\shorttitle{PLP Stabilization}
 
\date{\today}

\author[1]{Deepika Garg}
 \ead{deepika.lpu.pbi@gmail.com, deepikagarg@iisc.ac.in}

\author[1,2]{Sashikumaar Ganesan} 
  \ead{sashi@iisc.ac.in}
  
 \address[1]{National Mathematics Initiative, Indian Institute of Science, Bangalore - 560012, India}  
 \address[2]{Department of Computational and Data Sciences, Indian Institute of Science, Bangalore - 560012, India}
\begin{frontmatter}
\begin{abstract}
{\it A~priori} analysis for a generalized local projection stabilized conforming finite element approximation of  Darcy flow and Stokes problems is presented in this paper.  A first-order conforming $\mathbf{P}^{c}_{1}$  finite element space is used to approximate both the velocity and the pressure. It is shown that the stabilized discrete bilinear form satisfy the inf-sup condition with respect to a generalized local projection norm. Moreover, {\it a~priori} error estimates are derived for both   problems. Finally, the validation of the proposed stabilization scheme is demonstrated with appropriate numerical examples.
\end{abstract}
\end{frontmatter}
%%%%%%%%%%%%%%%%%%%%%%%%%%%%%%%%%%%%555555


%%%%%%%%%%%%%%%%%%%%%%%%%%%%%%%%%%%%%%%%%%%
\section{Introduction}
The numerical solution of Darcy equations has considerable practical importance in civil, petroleum, and electrical engineering, such as flow in porous media, heat transfer, semiconductor devices, etc. In general, numerical schemes for Darcy equations can be divided into two categories: (i) primal, a single-field formulation for pressure, and (ii) mixed two-field formulation in which pressure and velocity are variables.

Eliminate the velocity from mixed two-field formulation results in a scalar second-order partial differential equation ( PDEs) for the pressure. The construction of finite element methods based on this kind of formulation is straightforward. However, this direct approach results in lower-order velocity approximations compared to the pressure. Alternative approaches  such as mixed methods \cite{brezzi:2012:Fortin} and post-processing techniques \cite{loula:1995:Fernando}  have been used to improve the approximation of the velocity. The mixed finite element method based on the Galerkin formulation has    increasingly become popular   to discretize the Darcy   equations. The classical mixed variational formulation of Darcy equations is posed in the Sobolev spaces $\rm{H}(\mbox{div},\Omega) \ \text{and} \ {\rm{L}^{2}_{0}(\Omega)}$ for the velocity and pressure, respectively. It has been a challenge to develop  finite dimensional subspaces of these spaces that satisfy the {\it inf-sup} stability condition. Indeed, the choice of interpolation spaces is restricted when imposing this   {\it inf-sup} stability  condition. Nevertheless, a few finite element pairs that satisfy the {\it inf-sup} condition has been proposed. The well-known successful combinations are the Raviart-Thomas~\cite{Raviart:1977:Jean-Marie} and the Brezzi-Douglas-Marini \cite{BREZZI:MARINI}, which requires the continuity of normal component of the velocity in combination with specific discontinuous pressure interpolation. However,  such choices result in saddle point problems, which are more challenging to solve.

In this study, we propose  a mixed finite element formulation  with a generalized local projection stabilized {conforming finite element method} for Darcy equations, which avoids ${\rm H}(\mbox{div},\Omega)$ approximation space. 
% {This approach significantly simplifies the problem.}  
It is well-known that the application of standard Galerkin finite element method (FEM) to the Darcy equations induces spurious oscillations in the numerical solution. Nevertheless, the stability and accuracy of the standard Galerkin solution can be enhanced by applying a stabilization technique. Several stabilization methods such as streamline diffusion methods \cite{Tobiska:1990:sdfem,Tobiska:1996:Streamline}, least-square methods \cite{Brezzi:1993:bubble,bramble1971least,hughes1989new}, residual-free bubbles \cite{Brezzi:1993:bubble,Brezzi:1999:bubble,volker:2007:oseen}, local projection schemes \cite{braack:2006:MS,Dallmann:2016:LP,Stoke:2008:LP,Ganesan:2010:LPS,Nafa:2010:Localprojection}, continuous interior penalty methods \cite{Burman:2005:IP,BurmanErn:2007:CIPhp2,BurmanErn:2007:CIPhp1,Bur:2006:CIP} and many more  have been proposed in the literature. The basic idea of stabilization is to stabilize the Galerkin variational formulation so that the discrete approximation is stable and convergent; see, for example, \cite{Bochev:2008:Gunzburger,Bochev:2006,BurmanErn:2007:CIPhp2,Burman:2006:Hansbo2,Burman:2007:Hansbo1,Masud:2002:Thomas}. Stabilization methods for Stokes-like operators are well-studied in the literature, see for example~\cite{brack:2001:lp,Stoke:2008:LP} and a few studies for Darcy equations have also been presented, see for example \cite{Bochev:2008:Gunzburger,Bochev:2006,Masud:2002:Thomas,Nafa:2010:Localprojection}.

The local projection stabilization (LPS) method has been proposed in \cite{brack:2001:lp,braack:2006:MS} for the Stokes problem and subsequently extended to various other classes of problems ~\cite{Braack:2002:Schieweck,Stoke:2008:LP,Ganesan:2010:LPS,Guermond:2001,Knobloch:2010:LPoceen,Nafa:2010:Localprojection}. The LPS is based on a projection of the finite element space $Y_h$, which approximates the unknown to the discontinuous space $D_h$, see \cite{brack:2001:lp,braack:2006:MS}. LPS is very attractive, mainly because of its commutation properties in optimization problems \cite{Braack:2009:Optimal} and similar stabilization properties to those of residual approaches \cite{Knobloch:2011:Tobiska}. A significant benefit of the local projection method is that the LPS approach uses a symmetric stabilization term and contains fewer stabilization terms than the residual-based stabilization approach. Generalized local projection stabilization ({\rm GLPS}) is a more generalized form of LPS that allows us to define local projection spaces on overlapping sets. GLPS has first been introduced and studied for the convection-diffusion problem in \cite{ADTG,Knobloch:2010:LP} and for the Oseen problem in \cite{Biswas:Gudi,Knobloch:2010:LPoceen}, recently, for the advection-reaction equations in \cite{Ganesan:Deepika}. {\it A priori} analysis in \cite{Knobloch:2010:LP,Knobloch:2010:LPoceen} is based on an inf-sup condition for $Y_h$ and $D_h$ spaces and the existence of orthogonal projection of  $Y_h$ into $D_h$. Further, unlike LPS, GLPS needs neither a macro grid nor an enrichment of approximation spaces.

{\it The main contributions of this paper are the development of a GLPS conforming finite element scheme for Darcy equations and the derivation of its stability and convergence estimates.}  In the present analysis, we approximate the velocity and pressure with the piecewise linear polynomial finite element space. In particular, the use of piecewise linear finite elements for both the velocities and the pressure results in ill-posed discretizations. Therefore, GLPS is proposed in this work to suppress the oscillations in the approximations. The boundary conditions are not used strongly in discrete space; hence, the discrete formulation is a combination of standard Galerkin formulations, stabilization terms, and weakly imposed boundary conditions.  The proposed bilinear form satisfies an inf-sup condition with respect to generalized local projection stabilized norm, which leads to the well-posedness of the discrete problem.  {\it A priori} error analysis assures the optimal order of convergence, that is, $\mathcal{O}(h^{3/2})$ in the case of $(\mathbf{P}^{c}_1/\mathbf{P}^{c}_1)$ conforming finite element approximation. Furthermore, the above approach has also been used to study the Stokes problem. We give an elementary proof of stability and convergence analysis for the Stokes problem.

The outline of the article is as follows: In Section \ref{sec2}, we introduce the weak formulation of
the Darcy flow, notations, and preliminaries, which are used throughout the paper. Section \ref{sec3}
is devoted to an overlapping local projection stabilized conforming finite element methods in which we derive the
stability analysis with respect to a generalized local projection norm. In section \ref{sec4}, we provide an optimal {\it a priori} error estimates with respect to a generalized local projection norm. In Section {\ref{thestoke}}, we extended the above result to Stokes problem in the conforming FEM. Section  \ref{computation} presents some numerical experiments that confirm the theoretical analysis.

\section{{The Darcy problem}}\label{sec2}
Let $\Omega\subset\mathbb{R}^2$ be an open bounded polygonal domain with smooth boundary $\partial{\Omega}$. Consider the following Darcy flow equations: Find $(\textbf{u},p)$ such that
\begin{align}
\textbf{u} + \nabla p = \textbf{f}; \quad  \nabla \cdot\textbf{u} &= \phi \  \  \text{\ in} \  \Omega, \label{problem1} 
\\   
%\nabla \cdot\textbf{u} &= \phi  \ \  \text{\ in} \  \Omega,   \\
\textbf{u}\cdot \textbf{n} &= 0  \   \ \text{ on} \ \partial{\Omega} \nonumber. 
\end{align}
Here, $ \textbf{u} $ denotes the velocity vector, $p$ is the pressure, $ \textbf{f} \in [\rm{L}^{2}(\Omega)]^{2}$ is the source function, $\phi$ is the volumetric flow rate source,  and $\textbf{n}$ is the unit outward normal vector to $\partial{\Omega}$. The divergence constraint implies that the prescribed data must satisfy the condition %$\int_{\omega}\phi \dx=\int_{\Gamma}\psi \ds$ \\
\begin{align*}
\int_{\Omega}\phi \dx=0.
\end{align*}
In order to formulate a weak formulation of the Darcy flow equations, we consider the following Sobolev spaces
\[
\textbf{V} := \left\{\textbf{v} \in \rm{H}(\mbox{div},\Omega) |  \ \textbf{v} \cdot \textbf{n} =0 \ \text{on}  \ \partial{\Omega} \right\},  \quad 
Q := \rm{L}^{2}_{0}(\Omega)= \left\{q \in \rm{L}^{2}(\Omega) | \ \int_{\Omega}{q}\dx = 0 \  \right\},
\]
where $\rm{L}^{2}(\Omega)$ is a space of square-integrable measurable function. Moreover, a weak formulation of the model problem (\ref{problem1}) reads:
Find $(\textbf{u}, p)$ $\in \textbf{V} \times Q$ such that 
\[
a(\textbf{u},\textbf{v}) -b(p,\textbf{v}) = (\textbf{f},\textbf{v}); \qquad
b(\textbf{u},q)= (\phi,q),
\]
for all $v\in \textbf{V} $ and $q \in Q$. Here, $(\cdot,\cdot)$ denotes the $\rm{L}^{2}(\Omega)$ inner product and 
\begin{align*}
a(\textbf{u},\textbf{v}) := \int_{\Omega}{\textbf{u} \cdot \textbf{v}} \dx; \  \ \ b(p,\textbf{v}) := \int_{\Omega}p\nabla \cdot \textbf{v} \dx.
\end{align*}
An equivalent weak formulation of the model problem can be defined on the product space $ \textbf{V} \times Q $ and it reads:
Find $(\textbf{u}, p)   \ \in \textbf{V} \times Q \text{  such  that  }$
\begin{equation}\label{weak_f_d}
A((\textbf{u},p),(\textbf{v},q)) = L(\textbf{v}),
\end{equation}
 for  all  $(\textbf{v}, q) \in \textbf{V} \times Q$, where %\ \  $A((\textbf{u},p),(\textbf{v},q))= a(\textbf{u},\textbf{v}) -b(p,\textbf{v})+b(q,\textbf{u})$.\\
\[
A((\textbf{u},p),(\textbf{v},q)):= a(\textbf{u},\textbf{v}) -b(p,\textbf{v})+b(q,\textbf{u});\qquad
L(\textbf{v}):=(\textbf{f},\textbf{v})+(\phi,q).
\]
Furthermore, Banach-Ne$\check{c}$as-Babu$\check{s}$ka theorem \cite[pp. 85]{Ern:2004:FEM} guarantees that the model problem (\ref{problem1}) is well-posed in $\textbf{V} \times Q$, for more details; see \cite[pp. 230]{Ern:2004:FEM}.
%%%%%%%%%%%%%%%%%%%%%%%%%%%%%%%%%%%%%%%%%%%%%%%%%%%%%%%%
\subsection{Finite element formulation}
Let $\cT_{h}$ be a collection of non-overlapping quasi-uniform triangles obtained by a decomposition of $\Omega$.
Let $h_K=\rm{diam}(K)$  for all $K \in \cT_{h} $ and the mesh-size $ h = \mbox{max}_{K\in \cT_h} h_K$.  
Let $\cE_{h}=\cE_{h}^{I}\cup \cE_{h}^{B}$ be the set of all edges in  $\cT_{h}$, where  $\cE_{h}^{I}$ and $\cE_{h}^{B}$ are the set of all interior and boundary edges, respectively, and $h_E=\rm{diam}(E)$ for all $E\in \cE_{h}$. 
 Let $\cV_h:=\cV_{h}^I\cup\cV_{h}^B$ be the set of all vertices in  $\cV_{h}$, where  $\cV_{h}^{I}$ and $\cV_{h}^{B}$ are the set of all interior and boundary vertices, respectively. For any $a\in \cV_h$, we denote by ${\cM}_a$ (patch of $a$) the union of all cells that share the vertex $a$.
Further, define $h_a=\rm{diam}(\cM_{a})$ for all $a\in \cV_h$.
Moreover, {We use the following norm in the analysis. Let the piecewise constant function $h_{\cT}$ is defined by $h_{\cT} |_K = h_K$ and $s \in \mathbb{R}$ and $k \geq 0$}
\begin{align*}
\norm{h_{\cT}^{s} u}_{k} = \left(\sum_{K \in \cT_h } h_{K}^{2s} \norm{u}^{2}_{\rm{H}^{k}(K)} \right)^{\frac{1}{2}} \text{ for all} \ u \in \ \rm{H}^{k}(\cT_h).
\end{align*}
\begin{figure}[!ht]
\centering
\includegraphics[scale=0.5]{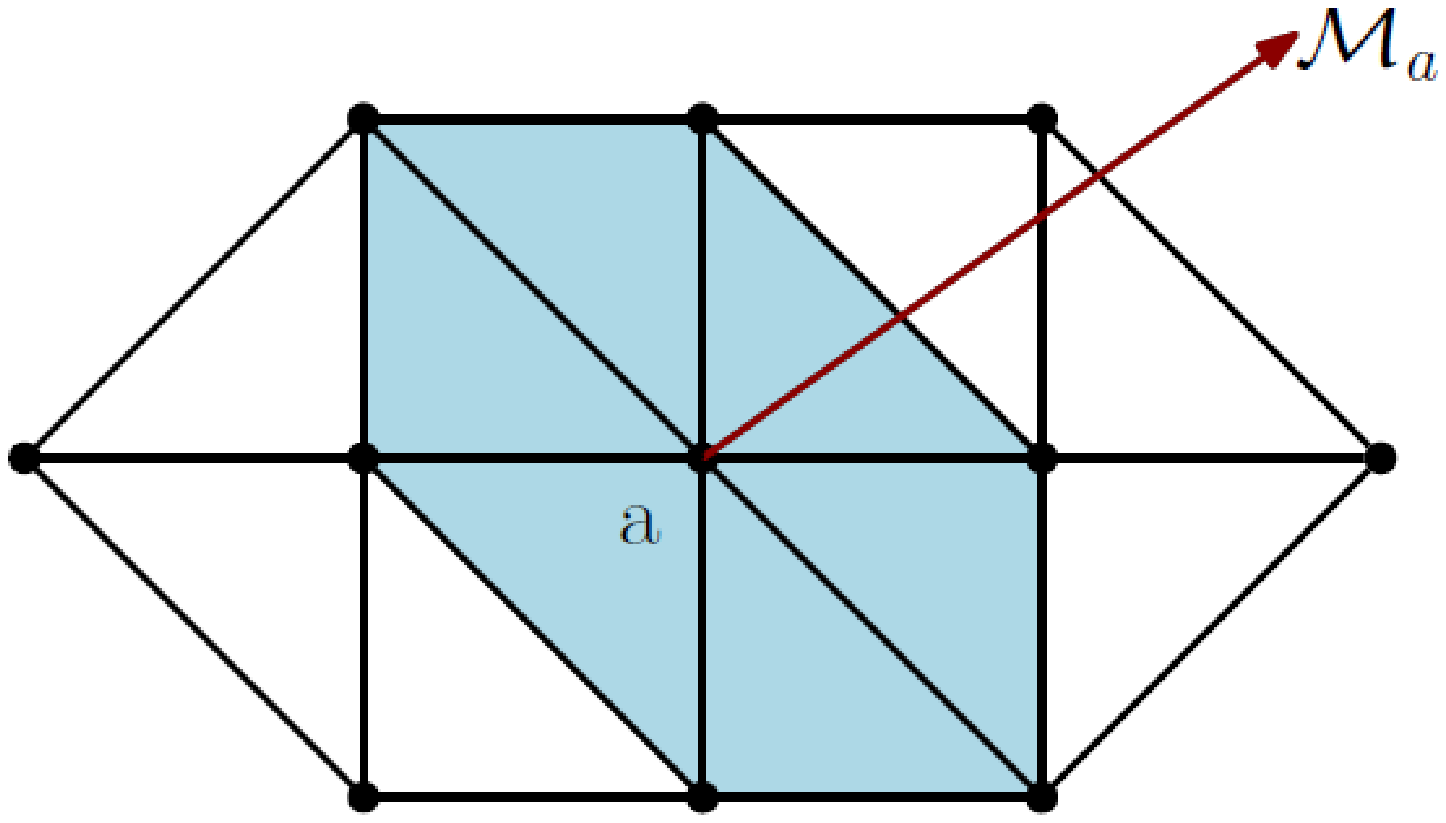}
 \caption{Node patch ${\cM}_a$.}
\end{figure}
{Suppose $I(a)$ denotes the index set for all $K_l$ elements, so that $K_l \subset {\cM} a$. } Then, the local mesh-size associated to ${\cM}_a$ is defined as
\begin{align*}
 {h}_a := \frac{1}{\mbox{card}(I(a))} \sum_{l\in I(a)} h_l, \quad \text{ for each}   \  a\in \cV_h,
\end{align*}
where $\text{card}(I(a))$ denotes the number of elements in ${\cM}_a$.
Since the mesh  $\cT_h$  is assumed to be locally quasi-uniform \cite{Bramble:2002:H1stb}, there exists a positive $\zeta \geq 1$ independent of $h$ such that
\begin{align*}
\zeta^{-1 }\leq \frac{ {h}_a}{ h_l} \leq \zeta    \text{ for all} \ l \in I(a).
\end{align*}
For any $a \in \cV_{h}$, define the fluctuation operator $\kappa_{a}: \rm{L}^{2}(\cM_{a}) \rightarrow \rm{L}^{2}(\cM_{a})$ by
\begin{align*}
\kappa_{a}(v)=v-\frac{1}{|\cM_{a}|} \int_{\cM_{a}}v \dx.
\end{align*}
%%%%%%%%%%%%%%%%%%%%%%%%%%%%%%%%
%%%%%%%%%%%%%%%%%%%%%%%%%%%%%%%%
 We next define a piecewise polynomial space as
\begin{align*}
\mathbb{P}_{k}(\cT_h):=\left\{v\in \rm{L}^{2}(\Omega): v|_K\in \mathbb{P}_k(K)\quad \forall K\in \cT_h\right\},
\end{align*}
where $\mathbb{P}_k(K)$, $k\ge 0$, is the space of polynomials of degree at most $k$ over the element $K$. Further, define a conforming finite element space of piecewise linear  
\begin{align*}
\mathbf{P}^{c}_{1}(\cT_h) := \left\{v \in \rm{H}^{1}(\Omega) \ : \ v|_K \in \mathbb{P}_1(K) ~~ \forall ~ K \in \cT_{h}   \right\}.
\end{align*}
Now recall the following technical results of finite element analysis.

\begin{lemma}{\rm Trace inequality} \cite[pp. 27]{Ern:2012:DGBook}: {Suppose E denotes an edge of $K \in \cT_{h}$}. For  $v_h\in \mathbb{P}_{k}(\cT_{h})$, there holds
	\begin{align} %\label{trace_ineq}
%	\|v\|_{\rm{L}^2(E)} &\leq C (h_K^{-1/2} \|v\|_{\rm{L}^2(K)} + h_K^{1/2} \|\nabla v\|_{\rm{L}^2(K)}), \\
	\|v_h\|_{\rm{L}^2(E)} &\leq C h_K^{-1/2} \|v_h\|_{\rm{L}^2(K)} . \label{trace_ineq1}
	\end{align}
\end{lemma}
\begin{lemma} {\rm Inverse inequality} \cite[pp. 26]{Ern:2012:DGBook}: Let $v \in \mathbb{P}_{k}(\cT_{h})$, for all $k \geq 0$; then
\begin{align}
\norm{\nabla{v}}_{K} \leq C h^{-1}_{K} \norm{v}_{K}.  \label{inverse_ineq1}
\end{align}
\end{lemma}

\begin{lemma}{\rm Poincar\'{e} inequality} \cite[pp. 104]{BScott:2008:FEM}: {For a bounded and connected  polygonal domain ${\Omega}$ and for any $v\in \rm{H}^1({\Omega})$, we have}
\begin{align*}%\label{poin} 
	\norm{v-\frac{1}{|{\Omega}|}\int_{{\Omega}}  v\dx}_{\rm{L}^2({\Omega})}\leq C h_{\Omega}\norm{\nabla v}_{\rm{L}^2({\Omega})},
	\end{align*}
	  {where $h_{\Omega}$ and $|{\Omega}|$ denote the diameter and the measure of domain ${\Omega}$}.
In particular, for every vertex $a\in \cV_h$ and every function $v \in \rm{H}^1({\cM_a})$, it holds
	\begin{align}\label{poin_2}
	\norm{v-\frac{1}{|{\cM}_a|}\int_{{\cM}_a}  v\dx}_{\rm{L}^{2}({\cM}_a)}\leq Ch_a\norm{\nabla v}_{\rm{L}^{2}({\cM}_a)},
	\end{align}
	where  the constant $C$ is independent of the mesh-size $h_a$.
\end{lemma}

%Let k $\geq $ 0 be an integer and set
{Furthermore, for a locally quasi-uniform and shape-regular triangulation the $\rm{L}^2$-orthogonal projection $I_h : \rm{L}^2 (\Omega) \rightarrow \mathbf{P}_{1}^{c}(\cT_{h})$ satisfies the following approximation properties}, for more details; see \cite{Yserentant:2014,ADTG}.

\begin{lemma}{$\rm{L}^2$-Orthogonal projections}: The $\rm{L}^2$-projection $I_h:\rm{L}^2(\Omega)\rightarrow  \mathbf{P}^{c}_{1}(\cT_h)$ satisfies
	\begin{align}\label{eq_30}
	\norm{ h_{\T}^{-1} ({v-I_{h}v})}+\norm{\nabla(v-I_hv)} &\leq C\norm{h_{\T}v}_{2},\ {\rm for~ all } \ v\in \rm{H}^2(\Omega),
	\end{align}
	For vector valued functions  $\mathbf{I}_{h}:[\rm{L}^2(\Omega)]^{2}\rightarrow  [\mathbf{P}^{c}_{1}(\cT_h)]^{2}$ satisfy  
	\begin{align}\label{intglobal}
	\norm{ h_{\T}^{-1} ({\mathbf{v}-\mathbf{I}_h\mathbf{v}})}+\norm{\nabla(\mathbf{v}-\mathbf{I}_h\mathbf{v})} &\leq C\norm{h_{\T}\mathbf{v}}_{2}\ {\rm for~ all }~ \mathbf{v}\in [\rm{H}^2(\Omega)]^{2}.
	\end{align}
	Moreover, the trace inequality over each edge imply
	\begin{align}\label{intglobal_trace}
	\left(\sum_{E\in\cE_h}\norm{\mathbf{v}-\mathbf{I}_h\mathbf{v}}^2_{\rm{L}^2(E)} \right)^{1/2}&\leq C\norm{h_{\T}^{3/2}\mathbf{v}}_{2}\ {\rm for~ all }~~ \mathbf{v}\in [\rm{H}^2(\Omega)]^{2}.
	\end{align}
	The orthogonality relation for all $\mathbf{v}_h\in [\mathbf{P}^{c}_{1}(\cT_h)]^{2}$ imply
	\begin{align}\label{eq:ortho}
	(\mathbf{v}-\mathbf{I}_h \mathbf{v},\mathbf{v}_h)_{\rm{L}^{2}(\Omega)}=0, 
	\end{align}
	The following approximation estimates hold for the $L^{2}$-orthogonal projection operator
%For the $L^{2}$-orthogonal projection operator, there holds
\begin{equation} \label{stab}
\norm{\mathbf{I}_h\mathbf{v}} \leq \norm{\mathbf{v}}, \quad 
 \norm{h^{-1}_{\cK}\mathbf{I}_h\mathbf{v}} \leq C \norm{h^{-1}_{\cK}\mathbf{v}}, \quad 
\norm{\nabla(\mathbf{I}_h{v})} \leq C\norm{\nabla{\mathbf{v}}}.
\end{equation}
\end{lemma}

Note that throughout this paper, C (sometimes subscripted) denotes a generic positive constant, which may depend on the shape-regularity of the triangulation but is independent of the mesh-size. Further, the notation $c \lesssim d $ represents the inequality $c \leq Cd$. 
 Moreover, $(\cdot,\cdot)$ represents the $\rm{L}^{2}(\Omega)$ inner product; and $\rm{L}^2(\Omega)$ and $\rm{L}^{\infty}(\Omega)$ norms are respectively denoted by $ \norm{u}$ and $ \norm{u}_{\infty} $.  The standard notation of Sobolev space $\rm{H}^{s}(\Omega),$ for s=1,2 and its norm $\norm{\cdot}_{r}$ respectively, are used. The notation $[\rm{L}^{2}(\Omega)]^{2}$ and $[\rm{H}^{1}(\Omega)]^{2},$ respectively, abbreviates the vector-valued version of $\rm{L}^{2}(\Omega)$ and $\rm{H}^{1}(\Omega)$ and $\rm{H}_{0}^{1}(\Omega)$ is a subspace of $\rm{H}^{1}(\Omega)$ with zero trace functions.

%%%%%%%%%%%%%%%%%%%%%%%%%%%%%%%%%%%%%%%%%%%%%%%%%%%%%%%%%%
\section{An overlapping local projection stabilization for Darcy flow problem }
This section describes the overlapping local projection stabilization conforming finite element methods for the problem (\ref{problem1}), where the velocity and the pressure are approximated with the continuous piecewise linear finite element spaces.
The velocity field will be sought in $ \textbf{V}_{h}:=[\mathbf{P}^{c}_{1}(\cT_h)]^{2}$ and the pressure in $Q_{h}:=\rm{L}_{0}^{2}(\Omega)\bigcap \mathbf{P}^{c}_{1}(\cT_{h})$.
An overlapping local projection stabilized conforming finite element method is defined as follows: 
Find $(\textbf{u}_h, p_h)\in \textbf{V}_{h} \times Q_h $ such that 
\begin{equation}\label{weak_f_cond} 
A_{h}((\textbf{u}_h,p_h),(\textbf{v},q)) = L(\textbf{v}, q), \ \text{for all}  \ (\textbf{v},q) \in \textbf{V}_{h} \times Q_h,
\end{equation} where
\begin{align}\label{bilinear_d}
  \ A_{h}((\textbf{u}_h,p_h),(\textbf{v},q)) = a_{h}(\textbf{u}_h,\textbf{v}) -b_{h}(p_h,\textbf{v})+b_{h}(\textbf{u}_h,q) + S_{h}((\textbf{u}_h,p_h),(\textbf{v},q)) ,
\end{align}  
and 
\begin{align*}
a_h(\textbf{u}_h,\textbf{v}) :&= \sum_{K \in \cT_{h}} \int_{K}{\textbf{u}_h\cdot\textbf{v}} \dx ,  \\ b_h(p_h,\textbf{v}) :&= (p_h,\nabla\cdot\textbf{v})-\sum_{E \in \cE_{h}^{B}}\int_{E}(\textbf{v} \cdot \textbf{n}) p_h \ds,  \\
S_{h}((\textbf{u}_h,p_h),(\textbf{v},q)):&=S_{si}((\textbf{u}_h,p_h),(\textbf{v},q))+S_{sb}((\textbf{u}_h,p_h),(\textbf{v},q)),  \\ 
S_{si}((\textbf{u}_h,p_h),(\textbf{v},q)) :&= \sum_{a \in \cV_{h}} \beta_{a} \int_{\cM_{a}}{\kappa_{a}(\nabla \cdot {\textbf{u}_h})}{\kappa_{a}(\nabla \cdot {\textbf{v}})} \dx \\ & \quad+\sum_{a \in \cV_{h}} \beta_{a} \int_{\cM_{a}}{\kappa_{a}(\nabla{p_h})}{\kappa_{a}(\nabla{q})} \dx,  \\
S_{sb}((\textbf{u}_h,p_h),(\textbf{v},q)) :&= \sum_{E \in \cE^{B}_{h}} \int_{E} ({\textbf{u}_h \cdot \textbf{n}})({\textbf{v} \cdot \textbf{n}}) \ds,\\
L(\textbf{v}, q) :&= (\textbf{f},\textbf{v}) + (\phi,q).
\end{align*}
Further, introduce a generalized local projection norm on $\textbf{V}_{h} \times Q_h $  by
\begin{equation}\label{GLP_norm_d}
\tnorm {(\textbf{u}_h,p_h)}^{2} := \Vert \textbf{u}_h \Vert^{2}+\Vert h^{\frac{1}{2}}_{\T} (\nabla \cdot \textbf{u}_h) \Vert^{2} + \Vert p_h \Vert^{2} +S_{h}( {(\textbf{u}_h,p_h)}, {(\textbf{u}_h,p_h)}).
\end{equation}
\subsection{{The inf-sup condition}} \label{sec3}
The main result of this section is the following theorem, which ensures that the discrete bilinear form is well-posed  \cite[pp. 85]{Ern:2004:FEM}.
\begin{theorem} \label{inf_sup_d_1} {The discrete bilinear form (\ref{weak_f_cond}) satisfies the following inf-sup condition for some positive constant $\gamma$, independent of $h$,  }  
\begin{equation*}
 \inf_{(\mathbf{u}_h,p_h) \in \textbf{V}_{h} \times Q_h } \sup_{(\mathbf{v}_h,q_h) \in \textbf{V}_{h} \times Q_h }\frac{A_{h}((\mathbf{u}_h,p_h),(\mathbf{v}_h,q_h))}{ \tnorm{(\mathbf{u}_h,p_h)}  \tnorm{ (\mathbf{v}_h,q_h) }} \geq \gamma.
\end{equation*}
%where the constant $\gamma$ is independent of the mesh-size $h$.
\end{theorem}
{\bf Proof.} In order to prove the stability result, it is enough to choose some $(\textbf{v}_h, q_h) \in \textbf{V}_{h} \times Q_h$ for any arbitrary $(\textbf{u}_h, p_h) \in \textbf{V}_{h} \times Q_h,$ such that
\begin{equation} \label{inf_sup_d} 
 \sup_{(\textbf{v}_h,q_h) \in \textbf{V}_{h} \times Q_h }\frac{A_{h}((\textbf{u}_h,p_h),(\textbf{v}_h,q_h))}{\tnorm {(\textbf{v}_h,q_h)}} \geq \gamma \tnorm{ (\textbf{u}_h,p_h)} > 0.
\end{equation}
We first consider the bilinear form in (\ref{bilinear_d}) with  $(\textbf{v}_h, q_h) = (\textbf{u}_h, p_h)$.  
\begin{equation} \label{eq_26a}
 {A_{h}((\textbf{u}_h,p_h),(\textbf{u}_h,p_h))} = \Vert \textbf{u}_h \Vert^{2}+ S_{h}((\textbf{u}_h,p_h),(\textbf{u}_h,p_h)).
\end{equation}
The stability of the pair $([{\rm{H}^{1}_{0}}(\Omega)]^2/\rm{L}_{0}^{2}(\Omega))$ \cite[pp. 199]{Ern:2004:FEM} implies that, there exists a constant $\mu>0$ such that
\begin{align} \label{dis_inf_sup_1}
\inf_{q_h \in  Q_h } \sup_{\mathbf{v} \in [{\rm{H}^{1}_{0}}(\Omega)]^2 }\frac{(\nabla\cdot{\textbf{v}}, q_h)}{\norm {\nabla\cdot{\textbf{v}}}\norm {q_h}} \geq \mu >0.
\end{align}
As a consequence of (\ref{dis_inf_sup_1}), for each $p_h \in Q_h$, there exists $\mathbf{z} \in [{\rm{H}^{1}_{0}}(\Omega)]^2$ such that
\begin{equation}\label{surjectivity_h} 
-(\nabla\cdot{\textbf{z}},p_h) = \norm{p_h}^{2} \  \text{and} \  \Vert \textbf{z} \Vert_{1,\Omega}  \leq C_1\Vert p_h \Vert.
\end{equation}
Let $\textbf{z}  \in  [\rm{H}_{0}^{1}(\Omega)]^{2}$ is defined as in (\ref{surjectivity_h}). Let $\textbf{z}_{h}=\textbf{I}_h\textbf{z} \in \textbf{V}_{h}$.  
\begin{align} \label{surjectivity_h1}
\norm{\textbf{z}_h}_{1, \Omega} \leq \norm{\textbf{z}}_{1, \Omega}\leq C_1\Vert p_h \Vert.
\end{align}
 Taking $(\textbf{v}_h, q_h) = (\textbf{z}_h, 0)$ as a test function pair, the bilinear form   
    (\ref{bilinear_d}) becomes
\begin{align}\label{z_h1}
A_{h}((\textbf{u}_h,p_h),(\textbf{z}_h, 0)) &= a_{h}(\textbf{u}_h,\textbf{z}_h)-b_{h}(p_h,\textbf{z}_h) + S_{h}((\textbf{u}_h,p_h),(\textbf{z}_h, 0)).
\end{align}
Let us now bound the three contributions. Applying Cauchy-Schwarz inequality, (\ref{surjectivity_h}) and the Young's inequality
\begin{align*}
\begin{split}
a_{h}(\textbf{u}_h,\textbf{z}_h) \leq \Vert \textbf{u}_h \Vert\Vert \textbf{z}_h \Vert
                \leq C_1 \Vert \textbf{u}_h \Vert\Vert p_h \Vert
                \leq C \Vert \textbf{u}_h \Vert^{2}+{\frac{1}{8}}\Vert p_h \Vert^{2}.
\end{split}
\end{align*}
In the second term of (\ref{z_h1}), add  $0=(p_h,p_h)-(p_h, -\nabla \cdot\textbf{z})$ to obtain
\begin{align}\label{22c}
-b_{h}(p_h,\textbf{z}_h)&=-(p_h,\nabla\cdot\textbf{z}_h)+\sum_{E \in \E_{h}^{B}}\int_{E}(\textbf{z}_h \cdot \textbf{n}) p_h \ds \nonumber \\ &= \norm{p_h}^{2} +(p_h,\nabla\cdot(\textbf{z}-\textbf{z}_h))+\sum_{E \in \cE_{h}^{B}}\int_{E}(\textbf{z}_h \cdot \textbf{n}) p_h \ds.
\end{align}
Applying an integration by parts to the second term of (\ref{22c}) we get
%Consider the second term of (\ref{22c}), an integration by parts leads to 
\begin{align*}% \label{eq_22a}
(p_h, \nabla\cdot({\textbf{z}-\textbf{z}_h})) = - (\nabla{p_h},({\textbf{z}-\textbf{z}_h})) +\sum_{E \in \cE_{h}^{B}}\int_{E}{p_h(\textbf{z}-\textbf{z}_h)\cdot \textbf{n}} \dx.
\end{align*}
It follows that
\begin{align*}%\label{22d}
-b_{h}(p_h,\textbf{z}_h)&= \norm{p_h}^{2} - (\nabla{p_h},{\textbf{z}-\textbf{z}_h}).
\end{align*}
Using the canonical nodal basis-function $\phi_a$ at the node $a \in \cV_{h}$ over the mesh ${\T}_h$. Since, $\sum_{a \in \cV_{h}}{\phi_{a}} \equiv 1 $, we have
%Let $\phi_a$ be the canonical nodel basis-function at the node $a \in \cV_{h}$ over the triangulation ${\T}_h$. Since $\sum_{a \in \cV_{h}}{\phi_{a}} \equiv 1 $ 
\begin{align} \label{basis_phi}
(\nabla{p_h},{\textbf{z}-\textbf{z}_h})   &= \sum_{K \in \cT_{h}}\int_{K} \nabla{p_h}({\textbf{z}-\textbf{z}_h})\sum_{a \in \cV_{h}}{\phi_{a}}\dx
%\\&= \sum_{a\in\cV_h}\sum_{K\in\cT_h} \int_{K} (\textbf{z}-\textbf{z}_h) \ \nabla p_h \phi_a \dx. \nonumber 
\\&= \sum_{a\in\cV_h} \int_{\cM_a} (\textbf{z}-\textbf{z}_h) \ \nabla p_h \phi_a \dx. \nonumber 
\end{align}
Using the orthogonality property of $\rm{L}^2$-projection (\ref{eq:ortho}) with the test function $C_a \phi_a \in \textbf{V}_h$ , where
$C_a=\frac{1}{|\cM_a|}\int_{\cM_a} \nabla p_h \dx$, and $\norm \phi_{\infty} \leq 1 $, we have
\begin{align*}
(\nabla{p_h},{\textbf{z}-\textbf{z}_h}) &=  \sum_{a\in\cV_h} \int_{\cM_a} (\textbf{z}-\textbf{z}_h) \Big( \nabla p_h - \frac{1}{|\cM_a|}\int_{\cM_a} \nabla p_h \dx \Big)\phi_a \dx \\ &\leq \left(\sum_{a\in\cV_h} \int_{\cM_a} \beta^{-1}_{a} (\textbf{z}-\textbf{z}_h)^{2} \dx\right)^{\frac{1}{2}}\left(\sum_{a\in\cV_h} \int_{\cM_a} \beta_{a} \kappa_{a}^{2}{( \nabla p_h)} \dx\right)^{\frac{1}{2}}
\\ \nonumber 
%&\leq \norm{h_{\T}\textbf{z}}_{2}[S^{nc}((\textbf{u}_h, p_h), (\textbf{u}_h, p_h))]^{\frac{1}{2}}\\
&\leq  \frac{1}{8}\norm{p_h}^{2}+ C S_{h}((\textbf{u}_h, p_h), (\textbf{u}_h, p_h)).
\end{align*}
Using the Cauchy-Schwarz inequality, the boundedness of an overlapping local projection operator and (\ref{surjectivity_h}) we obtain
\begin{align}\label{e_1}
S_{si}((\textbf{u}_h,p_h),(\textbf{z}_h,0)) &\leq [S_{si}((\textbf{u}_h,p_h),(\textbf{u}_h,p_h))]^{\frac{1}{2}}[S_{si}((\textbf{z}_h,0),(\textbf{z}_h,0))]^{\frac{1}{2}} \nonumber \\
                &\leq [S_{si}((\textbf{u}_h,p_h),(\textbf{u}_h,p_h))]^{\frac{1}{2}} \Vert \nabla\cdot\textbf{z}_h \Vert \nonumber \\
                &\leq  \frac{C}{2} S_{h}((\textbf{u}_h,p_h),(\textbf{u}_h,p_h)) + \frac{1}{8}\Vert p_h \Vert^{2}.                    
\end{align}
Since $\textbf{z}=0 $ on the boundary edges, using trace inequality over edges and (\ref{surjectivity_h}), the next term of stabilization is handled as
\begin{align*} 
S_{sb}((\textbf{z}_h,0),(\textbf{z}_h,0)) = \sum_{E \in \cE^{B}_{h}} \int_{E} ({\textbf{z}_h \cdot  \textbf{n}})^{2} \ds  &= \sum_{E \in \cE^{B}_{h}} \int_{E} ({(\textbf{z}-\textbf{z}_h) \cdot  \textbf{n}})^{2} \ds \nonumber \\ &\leq C \norm{h^{\frac{1}{2}}_{\T}\textbf{z}}_{1,\Omega} \leq \frac{1}{8} \norm{p_h}^{2}.
\end{align*}Thus,
\begin{align*} 
S_{h}((\textbf{u}_h,p_h),(\textbf{z}_h,0)) \leq \frac{C}{2} S_{h}((\textbf{u}_h,p_h),(\textbf{u}_h,p_h)) +\frac{1}{4} \norm{p_h}^{2}.
\end{align*} 
Put together, (\ref{z_h1}) leads to 
\begin{equation} \label{eq_25ab}
A_{h}((\textbf{u}_h,p_h),(\textbf{z}_h, 0)) \geq \frac{1}{2}\Vert p_h \Vert^{2} -C \Big(\Vert \textbf{u}_h \Vert^{2} + \frac{1}{2}S_{h}((\textbf{u}_h,p_h),(u_h,p_h)) \Big).
\end{equation}
Finally, the control of $\left\| {h^{\frac{1}{2}}_{\cT} (\nabla \cdot \textbf{u}_h)}\right\|^{2}$ can be obtained by choosing $(\textbf{v}_h,q_{h}) = (0,h_{\K} (\nabla \cdot \textbf{u}_h))$ in (\ref{bilinear_d}), that is,
\begin{align} \label{eq_23b}
  A_{h}((\textbf{u}_h,p_h),(0,I_{h}(h_{\T} (\nabla \cdot \textbf{u}_h)))) = b_{h}(&I_{h}(h_{\T} (\nabla \cdot \textbf{u}_h)),\textbf{u}_h) \nonumber\\&+S_{h}((\textbf{u}_h,p_h),(0,I_{h}(h_{\T} (\nabla \cdot \textbf{u}_h)))).
\end{align}
By adding and subtracting $\norm {h^{\frac{1}{2}}_{\T}(\nabla \cdot \textbf{u}_h)}^{2},$ the first term of (\ref{eq_23b}) becomes
\begin{align} \label{eq_23_a}
 b_{h}(I_{h}(h_{\T} (\nabla \cdot \textbf{u}_h)),\textbf{u}_h) &= \norm {h^{\frac{1}{2}}_{\T} (\nabla \cdot \textbf{u}_h)}^{2} + (I_{h}(h_{\T}(\nabla \cdot \textbf{u}_h)) - h_{\T}(\nabla \cdot \textbf{u}_h), \nabla \cdot \textbf{u}_h  ) \nonumber \\ & \quad-\sum_{E \in \cE_{h}^{B}}\int_{E}(\textbf{u}_h \cdot \textbf{n}) \ I_{h}(h_{\T} (\nabla \cdot \textbf{u}_h) )\ds.
\end{align} 
The second term of (\ref{eq_23_a}) is estimated as
\begin{align*}
(&I_h(h_{\T}(\nabla \cdot \textbf{u}_{h})) - h_{\T}(\nabla \cdot \textbf{u}_{h}), \nabla \cdot \textbf{u}_{h}  )
\\&= \sum_{a \in \cM_{a}}\int_{\cM_{a}} I_h(h_{K}(\nabla \cdot \textbf{u}_{h})) - h_{K}(\nabla \cdot \textbf{u}_{h})( \nabla \cdot \textbf{u}_{h} ) \phi_{a}\dx \\ &= \sum_{a \in \cM_{a}}\int_{\cM_{a}} \left(I_h(h_{K}(\nabla \cdot \textbf{u}_{h})) - h_{K}(\nabla \cdot \textbf{u}_{h})\right)\left( \nabla \cdot \textbf{u}_{h} -\frac{1}{|\cM_{a}|}\int_{\cM_{a}}\nabla \cdot \textbf{u}_{h}\dx \right) \phi_{a}\dx 
% \\&\leq \sum_{E \in \cE_{h}}\norm{ I_h(h_{\T}(\nabla \cdot \textbf{u}_{h})) - h_{\T}(\nabla \cdot \textbf{u}_{h})}\norm{( \nabla \cdot \textbf{u}_{h} )-\frac{1}{|\cM_{E}|}\int_{\cM_{E}}\nabla_{h} \cdot \textbf{u}_{h}\dx}  
\\ &\leq \left(\sum_{a \in \cM_{a}} \beta^{-1}_{a}\norm{ I_h(h_{\T}(\nabla \cdot \textbf{u}_{h})) - h_{\T}(\nabla \cdot \textbf{u}_{h})}_{\rm{L}^{2}(\cM_{a})}^{2}\right)^{\frac{1}{2}} [S_{si}((\textbf{u}_{h},0),(\textbf{u}_{h},0))]^{\frac{1}{2}} \\&\leq \frac{1}{6}\norm{h^{\frac{1}{2}}_{\T}( \nabla \cdot \textbf{u}_{h} )}^{2}+ \frac{C}{2} S_{si}((\textbf{u}_{h},p_h),(\textbf{u}_{h},q_h)).
\end{align*}
In the third term of (\ref{eq_23_a}), using the Cauchy-Schwarz inequality, the trace inequality, stability property of projection operator (\ref{stab}) and the Youngs inequality we get
\begin{align*}
\sum_{E \in \cE_{h}^{B}}\int_{E}(\textbf{u}_h \cdot \textbf{n}) \ & I_{h}(h_{K} (\nabla \cdot \textbf{u}_h))\ds \\ &\leq \left(\sum_{E \in \cE_{h}^{B}}\int_{E}(\textbf{u}_h \cdot \textbf{n})^{2}\ds\right)^{\frac{1}{2}} \left(\sum_{E \in \cE_{h}^{B}}\int_{E}(I_{h}(h_{K} (\nabla \cdot \textbf{u}_h))^{2}\ds\right)^{\frac{1}{2}} \\ &\leq   \frac{1}{6} \norm{h^{\frac{1}{2}}_{\T} (\nabla \cdot \textbf{u}_h)}^{2}+ \frac{C}{4} S_{sb}((\textbf{u}_h,0),(\textbf{u}_h,0)).
\end{align*}
Put together, (\ref{eq_23b}) leads to 
\begin{equation} \label{eq_24ab}
A_{h}((\textbf{u}_h,p_h),(0, I_h(h_{\T} (\nabla \cdot \textbf{u}_h)))) \geq \frac{1}{2}\norm {h^{\frac{1}{2}}_{\T} (\nabla \cdot \textbf{u}_h)}^{2} -\frac{C}{2}  S_{h}((\textbf{u}_h,p_h),(\textbf{u}_h,p_h)).
\end{equation}
The selection of $(\textbf{v}_h, q_h)$ is
\begin{align*}
(\textbf{v}_h, q_h) = (\textbf{u}_h, p_h) +\frac{1}{{C}+1}(\textbf{z}_h, 0)+\frac{1}{{C}+1}(0, I_h(h_{\T} \nabla \cdot \textbf{u}_h)).
\end{align*}
Adding (\ref{eq_26a}), (\ref{eq_25ab}) and (\ref{eq_24ab}) leads to 
\begin{align}\label{eq_27_a}
  A_{h}((&\textbf{u}_h,p_h),(\textbf{v}_h,q_h))\nonumber \\ &\geq  \Vert \textbf{u}_h \Vert^{2}+ S_h((\textbf{u}_h,p_h),(\textbf{u}_h,p_h))+\frac{1}{2+2{C}}\Vert p_h \Vert^{2} \nonumber \\ &\quad -\frac{{C}}{{C}+1} \left(\Vert \textbf{u}_h \Vert^{2} +\frac{1}{2}S_h((\textbf{u}_h,p_h),(\textbf{u}_h,p_h))\right) +\frac{1}{2+2{C}}\norm {h^{\frac{1}{2}}_{\T} (\nabla \cdot \textbf{u}_h)}^{2} \nonumber\\&\quad -\frac{{C}}{2+2{C}} \left(S_h((\textbf{u}_h,p_h),(\textbf{u}_h,p_h)) \right)\nonumber \\ 
  &=\frac{1}{2+2{C}}\Vert p_h \Vert^{2}+\frac{1}{2+2{C}}\norm {h^{\frac{1}{2}}_{\T} (\nabla \cdot \textbf{u}_h)}^{2} \nonumber\\&\quad +\left(1-\frac{{C}}{1+{C}}\right) \left(\Vert \textbf{u}_h \Vert^{2}+ S_h \big((\textbf{u}_h,p_h),(\textbf{u}_h,p_h) \big) \right) \nonumber \\
  &= \frac{1}{2+2{C}}\Vert p_h \Vert^{2}+\frac{1}{2+2{C}}\norm {h^{\frac{1}{2}}_{\T} (\nabla \cdot \textbf{u}_h)}^{2} \nonumber\\&\quad +\frac{1}{{C}+1} \Big(\Vert \textbf{u}_h \Vert^{2}+ S_h\left((\textbf{u}_h,p_h),(\textbf{u}_h,p_h)\big)\right) \nonumber \\
 &\geq \frac{1}{2{C}+2} \tnorm{(\textbf{u}_h,p_h)}^{2}.  
\end{align}
Applying the triangle inequality 
\begin{align} \label{eq_28a}
\tnorm{ (\textbf{v}_h,q_h) } &\leq \tnorm{ (\textbf{u}_h,p_h) }+\frac{1}{{C}+1}\tnorm{ ({\textbf{z}_h, 0}) }+\frac{1}{{C}+1}\tnorm{ ({0, h_{\T} (\nabla \cdot \textbf{u}_h)}) } \leq  {\alpha}\tnorm{ (\textbf{u}_h,p_h) }.
\end{align}
In the second term of (\ref{eq_28a}), applying  (\ref{surjectivity_h}) and a similar technique in (\ref{e_1}), we get
\begin{align*}
\tnorm{( {\textbf{z}_h, 0})} &= \Vert \textbf{z}_h \Vert^{2}+\Vert h^{\frac{1}{2}}_{\T} (\nabla \cdot \textbf{z}_h) \Vert^{2}  +S_{h}(\textbf{z}_h,\textbf{z}_h)\leq C\Vert {p}_h \Vert^{2}, 
\end{align*}
and in the third term of (\ref{eq_28a}), an inverse inequality (\ref{inverse_ineq1}) result in
\begin{align*}
\tnorm{ ({0, h_{\T} (\nabla \cdot \textbf{u}_h)}) } &= \Vert h_{\T} (\nabla \cdot \textbf{u}_h) \Vert^{2}   \leq C\Vert { \textbf{u}}_h \Vert^{2}. 
\end{align*}
Finally, (\ref{eq_27_a}) and (\ref{eq_28a}) lead to (\ref{inf_sup_d}), and these concludes the proof.

\subsection{{A priori error estimates}} \label{sec4}
This section presents {\it a priori} error estimates for the
$[\mathbf{P}^{c}_{1} /\mathbf{P}^{c}_{1}]$ approximation for velocity-pressure pair with respect to the $\tnorm{\cdot}$ norm.
\begin{lemma} \label{eq_01} Suppose $\beta_a = \beta{h_a};$ for some $\beta > 0.$ Let $(\mathbf{u}, p) \in [{\rm{H}^{2}}(\Omega)]^{2}\times  \rm{L}^{2}_{0} \bigcap \rm{H}^{2}(\Omega)$. Then
\begin{align} \label{eq_31}
\tnorm {({\mathbf{u}}-\mathbf{I}_h {\mathbf{u}}, p-I_h{p})} \leq C \left(\norm{h_{\T}^{\frac{3}{2}} \mathbf{u}}_{2} + \norm{h_{\T}^{\frac{3}{2}} p}_{2} \right ).
\end{align}
\end{lemma}
{\bf Proof.} We first consider the terms in $\tnorm{\cdot}$ norm defined in (\ref{GLP_norm_d}) 
\begin{align*}
\tnorm {(\textbf{u}-\mathbf{I}_h{\textbf{u}},p-I_hp)}^{2} := \Vert \textbf{u}-\textbf{I}_h{\textbf{u}} \Vert^{2}+\Vert h^{\frac{1}{2}}_{\K} (\nabla \cdot (\textbf{u}-\textbf{I}_h{\textbf{u}})) \Vert^{2} + \Vert p-I_h{p} \Vert^{2} \\ +S_{h}((\textbf{u}-\textbf{I}_h{\textbf{u}},p-I_hp),(\textbf{u}-\textbf{I}_h{\textbf{u}},p-I_hp)).
\end{align*}
Using  the projection estimates (\ref{eq_30})-(\ref{intglobal})  we get
\begin{align*}
\Vert \textbf{u}-\textbf{I}_h{\textbf{u}} \Vert \leq \norm{h_{\T}^{2} \textbf{u}}_{2}  , \  \Vert h^{\frac{1}{2}}_{\T}(\nabla\cdot(\textbf{u}-\textbf{I}_h{\textbf{u}})) \Vert \leq \norm{h^{\frac{3}{2}}_{\T} \textbf{u}}_{2} \ \text{and}  \  \norm {p-I_h p} \leq  \norm{h_{\T}^{2} p}_{2}. 
\end{align*}
Recall the stabilization term 
\begin{align} \label{new}
S_{h}((\textbf{u}-&\textbf{I}_h{\textbf{u}}, p-I_h{p}),(\textbf{u}-\textbf{I}_h{\textbf{u}}, p-I_h{p})) \nonumber \\
 &= \sum_{a \in \cV_{h}} \beta_{a} \int_{\cM_{a}}{\kappa^{2}_{a}(\nabla \cdot (\textbf{u}-\textbf{I}_h{\textbf{u}}))} \dx +\sum_{a \in \cV_{h}} \beta_{a} \int_{\cM_{a}}{\kappa^{2}_{a}(\nabla  (p-I_h{p}))} \dx \nonumber \\
&+ \sum_{E \in \cE^{B}_{h}} \int_{E} ({(\textbf{u}-\textbf{I}_h{\textbf{u}}) \cdot \textbf{n}})^{2} \ds.
\end{align}
In the first term of (\ref{new}), using the boundedness of an overlapping local projection operator, $\beta_a = \beta h_a$ and (\ref{intglobal}) we obtain
\begin{align*}
 \sum_{a \in \cV_{h}} \beta_{a} \int_{\cM_{a}}&{\kappa^{2}_{a}(\nabla \cdot (\textbf{u}-\textbf{I}_h{\textbf{u}}))} \dx 
 \\&= \sum_{a \in \cV_{h}} \beta_a \norm{\nabla \cdot (\textbf{u}-\textbf{I}_h{\textbf{u}})- 
\frac{1}{|\cM_a|} \int_{\cM_{a}}{\nabla\cdot(\textbf{u}-\textbf{I}_h{\textbf{u}})}\dx}_{\rm{L}^{2}(\cM_{a})}^{2} \nonumber \\   &\leq \sum_{a \in \cV_{h}} \beta_a  \norm { {\nabla \cdot(\textbf{u}-\textbf{I}_h{\textbf{u}})}}_{\rm{L}^{2}(\cM_{a})}^{2}    
 %\leq \sum_{z \in \cV_{h}} \delta h_{a} \norm{\nabla\cdot (\textbf{u}-\textbf{J}_h{\textbf{u}})}^{2}    
  \nonumber \\
 &\leq C \norm {h_{\T}^{3/2}\textbf{u}}^{2}_{2}
\end{align*} 
Similarly,
\begin{align*}
 \sum_{a \in \cV_{h}} \beta_{a} \int_{\cM_{a}}{\kappa^{2}_{a}(\nabla  (q-I_h{q}))} \dx
 \leq C \norm {h_{\T}^{3/2}q}^{2}_{2}
\end{align*} 
The boundary term is handled by using the trace inequality over each edges (\ref{intglobal_trace})
\begin{align*}
\sum_{E \in \cE^{B}_{h}} \int_{E} ({(\textbf{u}-\textbf{I}_h{\textbf{u}}) \cdot \textbf{n}})^{2} \ds=\sum_{E \in \cE^{B}_h} \norm{(\textbf{u}-\textbf{I}_h{\textbf{u}})\cdot \textbf{n}}_{\rm{L}^2(E)}^{2} \leq \norm {h_{\T}^{3/2} \textbf{u}}^{2}_{2}.
\end{align*}
The combination of above estimates leads to (\ref{eq_31}). This concludes the proof.
\begin{lemma} \label{eq+} Suppose $\beta_a = \beta{h_a};$ for some $\beta > 0.$ Let $(\mathbf{u}, p) \in [{\rm{H}^{2}}(\Omega)]^{2}\times  \rm{L}^{2}_{0} \bigcap \rm{H}^{2}(\Omega)$  and for all $(\mathbf{v}_h, q_h) \in \textbf{V}_{h} \times Q_h.$ Then
\begin{equation}\label{eq_35}
A_{h}((\mathbf{u}-\mathbf{I}_h{\mathbf{u}}, p-I_h{p}),(\mathbf{v}_h, q_h)) \leq C \left(\norm{h_{\T}^{\frac{3}{2}} \mathbf{u}}_{2} + \norm{h_{\T}^{\frac{3}{2}} p}_{2}\right ) \tnorm {(\mathbf{v}_h, q_h)}.
\end{equation}
\end{lemma}
{\bf Proof.} Consider the bilinear form in (\ref{bilinear_d}) 
\begin{align} \label{eq_32}
A_{h}((\textbf{u}-\textbf{I}_h{\textbf{u}}, p-I_h{p}),&(\textbf{v}_h, q_h))\nonumber \\ &= a_{h}(\textbf{u}-\textbf{I}_h{\mathbf{u}},\mathbf{v}_h) -b_{h}(p-I_h{p},\textbf{v}_h)+b_{h}(\textbf{u}-\textbf{I}_h{\textbf{u}},q_h)  \nonumber \\ &+ S_{h}((\textbf{u}-\textbf{I}_h{\textbf{u}}, p-I_h{p}),(\textbf{v}_h, q_h)).
\end{align}
Applying Cauchy-Schwarz inequality and the $\rm{L}^2$-projection property (\ref{intglobal}) 
\begin{align*}
a_{h}(\textbf{u}-\textbf{I}_h{\textbf{u}},\textbf{v}_h) \leq \norm{\textbf{u}-\textbf{I}_h{\textbf{u}}} \norm{\textbf{v}_h} \leq \norm{h_{\T}^{2} \textbf{u}}_{2} \tnorm{(\textbf{v}_h, q_h)}.
\end{align*}
Consider the second term of bilinear form (\ref{eq_32})
\begin{align} \label{eq_34}
b_{h}(p-I_h{p},\textbf{v}_h) &=  (p-I_h{p},\nabla \cdot\textbf{v}_h)-\sum_{E \in \E_{h}^{B}}\int_{E}(\textbf{v}_h \cdot \textbf{n}) \ (p-I_h{p})\ds.
\end{align}
Using Cauchy-Schwarz inequality and the $\rm{L}^2$-projection property in the first term of (\ref{eq_34}) we obtain
\begin{align*} %\label{eq_34b}
 (p-I_h{p},\nabla \cdot\textbf{v}_h) \leq \norm{p-I_h{p}} \norm{\nabla \cdot\textbf{v}_h} 
 \leq \norm{h_{\T}^{\frac{3}{2}}p} \norm{h_{\T}^{\frac{1}{2}}(\nabla \cdot\textbf{v}_h)}  
 \leq \norm{h_{\T}^{\frac{3}{2}}p} \tnorm{(\textbf{v}_h,q_h)}.
\end{align*}
 The second term of (\ref{eq_34}) is handled by using the Cauchy Schwarz inequality and trace inequality over edges, 
\begin{align*}
\sum_{E \in \cE_{h}^{B}}\int_{E}(\textbf{v}_h \cdot \textbf{n}) \ (p-I_h{p})\ds &\leq  \Big(\sum_{E \in \cE_{h}^{B}}\norm{\textbf{v}_h \cdot \textbf{n}}^{2}_{\rm{L}^{2}(E)}\Big)^{\frac{1}{2}} \Big(\sum_{E \in \cE_{h}^{B}} \norm{p-I_h{p}}^{2}_{\rm{L}^{2}(E)} \Big)^{\frac{1}{2}} \\ \nonumber &\leq \norm{h_{\T}^{\frac{3}{2}}p} \tnorm{(\textbf{v}_h,q_h)}.
\end{align*} 
Applying an integration by parts in the next term of the bilinear form (\ref{eq_32})
\begin{align}\label{eq_33}
b_{h}(\textbf{u}-\textbf{I}_h{\textbf{u}},q_h) 
&= (q_h, \nabla\cdot({\textbf{u}-\textbf{I}_h{\textbf{u}}}))-\sum_{E \in \cE_{h}^{B}}\int_{E}(({\textbf{u}-\textbf{I}_h{\textbf{u}}})\cdot \textbf{n} )\ q_h \ds \nonumber \\ 
&= -(\nabla{q_h},{\textbf{u}-\textbf{I}_h{\textbf{u}}} )+\sum_{E \in \cE_{h}^{B}}\int_{E}(({\textbf{u}-\textbf{I}_h{\textbf{u}}})\cdot \textbf{n} )\ q_h \ds \nonumber \\&-\sum_{E \in \cE_{h}^{B}}\int_{E}(({\textbf{u}-\textbf{I}_h{\textbf{u}}})\cdot \textbf{n} ) q_h \ds.
\end{align}
Applying the similar techniques as in (\ref{basis_phi}), the first term of (\ref{eq_33}) is estimated as:
\begin{align*} %\label{eq_34c}
(\nabla{q_h},{\textbf{u}-\textbf{I}_h{\textbf{u}}} )
\leq  C \norm{h^{\frac{3}{2}}_{\T} \textbf{u}}_{2} \tnorm{(\textbf{v}_h,q_h)}.
\end{align*}
%In the stabilization terms, again use of the Cauchy-Schwarz inequality and 
The last term is estimated in a similar way as in (\ref{new}) 
\begin{align*}
 S_{h}\big((\textbf{u}-\textbf{I}_h{\textbf{u}}, p-I_h{p}),(\textbf{v}_h, q_h)\big) 
 %&= F^{c/2}_{h}((\textbf{u}-\textbf{I}_h{\textbf{u}}, p-I_h{p}),(\textbf{u}-\textbf{I}_h{\textbf{u}}, p-%I_h{p}) )F^{c/2}_{h} ((\textbf{v}_h, q_h), (\textbf{v}_h, q_h)) \nonumber \\ 
 &\leq \Big(\norm{h_{\T}^{\frac{3}{2}} \textbf{u}}_{2} + \norm{h_{\T}^{\frac{3}{2}} p}_{2}\Big) \tnorm{(\textbf{v}_h, q_h)}.
\end{align*}
The collection of all above estimates shows (\ref{eq_35}) and concludes the proof.
\begin{lemma}{\rm Consistency Error}:\label{lm10} Suppose $(\mathbf{u}, p) \in [{\rm{H}^{2}}(\Omega)]^{2}\times  \rm{L}_{0}^{2}(\Omega) \bigcap \rm{H}^{2}(\Omega)\  and \\ \ (\mathbf{u}_h, {p}_h) \in \textbf{V}_h \times\  Q_h $ be the solutions to (\ref{weak_f_d}) and (\ref{weak_f_cond}), respectively. For any $(\mathbf{v}_h, q_h) \in \textbf{V}_h \times Q_h.$ Then
\begin{equation} \label{eq_36}
A_{h}((\mathbf{u}-\mathbf{u}_h, p-p_h),(\mathbf{v}_h, q_h)) \leq C \Big(\norm{h_{\T}^{\frac{3}{2}} \mathbf{u}}_{2} + \norm{ h_{\T}^{\frac{3}{2}} p}_{2} \Big) \tnorm {(\mathbf{v}_h, q_h)}.
\end{equation}
\end{lemma}
{\bf Proof.}
The model problem with the test function $(\textbf{v}_h, q_h) \in \textbf{V}_h \times Q_h$
and the definition of the bilinear form and the fact that the normal component $\textbf{u} \cdot \textbf{n}=0$ over the boundary edges
\begin{align*}
A_{h}((\textbf{u}-\textbf{u}_h, p-p_h),(\textbf{v}_h, q_h)) =S_{si}((\mathbf{u}, p),(\mathbf{v}_h, q_h)).
% -\sum_{e \in \xi^{b}_{h}} \int_{e}{\textbf{v}_h np}.
\end{align*}
\begin{align} \label{eq_37}
S_{si}((\mathbf{u}, p),(\mathbf{v}_h, q_h)) = \sum_{a \in \cV_{h}} \beta_{a} \int_{\cM_{a}}{\kappa_{a}(\nabla\cdot{\mathbf{u}})}&{\kappa_{a}(\nabla\cdot{\mathbf{v}_h})}\dx \nonumber \\ &+\sum_{a \in \cV_{h}} \beta_{a} \int_{\cM_{a}}{\kappa_{a}(\nabla p)} \cdot  {\kappa_{a}(\nabla {q_h})}\dx.
\end{align} 
Using the Cauchy-Schwarz inequality, the
Poincar$\acute{e}$ inequality (\ref{poin_2}) and $\beta_a = \beta h_a $  in the first term of (\ref{eq_37}) we have
\begin{align*}
\sum_{a \in \cV_{h}} \beta_{a} \int_{\cM_{a}}{\kappa_{a}(\nabla\cdot{\mathbf{u}})}&{\kappa_{a}(\nabla\cdot{\mathbf{v}_h})}\dx \nonumber \\&\leq \left(\sum_{a \in \cV_{h}} \beta_{a} \norm{\nabla \cdot \mathbf{u}- \frac{1}{|\cM_a|}\int_{\cM_a}{\nabla \cdot \mathbf{u}}\dx}_{\rm{L}^{2}(\cM_{a})}^{2}\right)^{1/2} 
 \nonumber \\ &S_{h}^{1/2}((\mathbf{u}, p),(\mathbf{v}_h, q_h))  \nonumber \\ &\leq \norm {h_{\T}^{3/2}\textbf{u}}_{2} \tnorm{(\mathbf{v}_h, q_h)}.
\end{align*}
In a similar way, the second  term is handled as:
\begin{align*} %\label{eq_38}
\sum_{a \in \cV_{h}} \beta_{a} \int_{\cM_{a}}{\kappa_{a}(\nabla p)}{\kappa_{a}(\nabla {q_h})} \dx\leq \norm {h_{\T}^{3/2}p}_{2} \tnorm{(\mathbf{v}_h, q_h)}.
\end{align*}
The collection of all above estimates shows (\ref{eq_36}) and concludes the proof. 
\begin{theorem} \label{th_01}
 Let $(\mathbf{u}, p) \in [{\rm{H}^{2}}(\Omega)]^{2}\times  \rm{L}^{2}_{0} \bigcap H^{2}(\Omega)\  and \ \ (\mathbf{u}_h, {p}_h) \in \textbf{V}_h \times\  Q_h$ be the solutions to (\ref{weak_f_d}) and (\ref{weak_f_cond}), respectively. Suppose $\beta_a = \beta{h_a};$ for some $\beta > 0.$ Then it holds  
  \begin{equation} \label{eq_40}
  \tnorm {(\mathbf{u}-\mathbf{u}_h, p-p_h)} \leq C \Big(\norm{h_{\T}^{\frac{3}{2}} \mathbf{u}}_{2} + \norm{ h_{\T}^{\frac{3}{2}} p}_{2} \Big).
  \end{equation}
  \end{theorem}
{ \bf Proof.} The triangle inequality implies
\begin{equation} \label{triangle}
\tnorm {(\mathbf{u}-\mathbf{u}_h, p-p_h)} \leq \tnorm {(\mathbf{u}- \textbf{I}_h {\mathbf{u}}, p-I_h{p})} + \tnorm {(\textbf{I}_h{\mathbf{u}}-\mathbf{u}_h,I_h{p}-p_h)}.
\end{equation}
The first term of (\ref{triangle}) follows from Lemma \ref{eq_01} i.e.
\begin{equation*}
\tnorm {(\textbf{u}-\textbf{I}_h{\textbf{u}}, p-I_h{p})} \leq C \Big(\norm{h_{\T}^{\frac{3}{2}} \textbf{u}}_{2} + \norm{ h_{\T}^{\frac{3}{2}} p}_{2}\Big).
\end{equation*}
The second of (\ref{triangle}) is handled by using Theorem \ref{inf_sup_d_1}
\begin{align} \label{eq_39}
\tnorm {(\textbf{I}_h{\textbf{u}}-\textbf{u}_h,I_h{p}-p_h)} &\leq {1/ \beta} \sup_{(\textbf{v}_h,q_h) \in \textbf{V}_h \times\  Q_h }\frac{A_{h}((\textbf{I}_h{\textbf{u}}-\textbf{u}_h,I_h{p}-p_h),(\textbf{v}_h,q_h))}{\tnorm{{\textbf{v}_h,q_h}}} \nonumber \\
&\leq {1/ \beta} \sup_{(\textbf{v}_h,q_h) \in \textbf{V}_h \times\  Q_h }\frac{A_{h}(\textbf{u}-\textbf{u}_h,p-p_h),(\textbf{v}_h,q_h))}{\tnorm{{\textbf{v}_h,q_h}}}  \nonumber \\
 &+\sup_{(\textbf{v}_h,q_h) \in \textbf{V}_h \times\  Q_h }\frac{A_{h}((\textbf{I}_h{\textbf{u}}-\textbf{u},I_h{p}-p),(\textbf{v}_h,q_h))}{\tnorm{{\textbf{v}_h,q_h}}}.
\end{align}
Finally, the result follows by using Lemma \ref{eq+} and Lemma \ref{lm10} in (\ref{eq_39}) and this concludes the proof.

\section{ The Stokes problem }\label{thestoke}
In this section, we extend the above analysis   to the Stokes  problem on mixed form. Consider the following Stokes  problem: 
\begin{align} 
-\Delta\textbf{u} + \nabla p  = \textbf{f}; \quad \nabla \cdot\textbf{u} &= 0   \  \  \text{\ in} \  \Omega, \label{stoke}\\
\textbf{u} &= 0  \   \ \text{ on} \ \partial{\Omega}. \nonumber 
\end{align}
where $\Omega\subset\mathbb{R}^2$ be an open bounded polygonal domain with smooth boundary $\partial{\Omega}$. Here, $ \textbf{u} $ denotes the velocities, $p$ denotes the pressure, $ \textbf{f} \in [\rm{L}^{2}(\Omega)]^{2}$ is some given data. The weak form of Stokes problem is obtained by considering the bilinear form
\begin{align*}
B((\textbf{u},p),(\textbf{v},q)):&= a(\textbf{u},\textbf{v}) -b(p,\textbf{v})+b(q,\textbf{u}),
\end{align*}
where $a(\textbf{u},\textbf{v})=(\nabla{\textbf{u}},\nabla{\textbf{v}})$ and $b(p,\textbf{v})=(p, \nabla \cdot \textbf{v}).$ We consider the functional spaces $\mathbf{V}=\{ \textbf{v} \in [\rm{H}^{1}_{0}(\Omega)]^{2}\}$, and $Q =\rm{H}^{1}(\Omega) \cap \rm{L}^{2}_{0}(\Omega)$. The weak formulation of (\ref{stoke}) now writes: Find $(\textbf{u},p) \in V \times Q$ such that
\begin{align*}
B((\textbf{u},p),(\textbf{v},q))=(\textbf{f},\textbf{v}).
\end{align*}
The existence of a weak solution to this problem follows by the application of the Lax-Milgram lemma in the divergence free subspace of $\mathbf{V}$ and the pressure in $Q$ by the Brezzi condition \cite{Brezzi:1984}.
%\subsection{An Overlapping Local Projection stabilization for Stokes problem }

Now, we describe an overlapping local projection stabilization conforming finite element methods for the problem (\ref{stoke}), where we approximate the velocity and the pressure with the continuous piecewise linear finite element spaces.
The velocity field will be sought in $ \textbf{V}_{h}:=[\mathbf{P}^{c}_{1}(\cT_h)]^{2}$ and the pressure in $Q_{h}:=\rm{L}_{0}^{2}(\Omega)\bigcap \mathbf{P}^{c}_{1}(\cT_{h})$.
An overlapping local projection stabilized conforming finite element method is defined as follows: 
Find $(\textbf{u}_h, p_h)\in \textbf{V}_{h} \times Q_h $ such that 
\begin{equation}\label{weak_f_cond_1} 
B_{h}((\textbf{u}_h,p_h),(\textbf{v},q)) = L(\textbf{v}, q), \ \text{for all}  \ (\textbf{v},q) \in \textbf{V}_{h} \times Q_h,
\end{equation} where
\begin{align}\label{bilinear_s}
   B_{h}((\textbf{u}_h,p_h),(\textbf{v},q)) = a_{h}(\textbf{u}_h,\textbf{v}) -b_{h}(p_h,\textbf{v})+b_{h}(\textbf{u}_h,q) + S_{h}((\textbf{u}_h,p_h),(\textbf{v},q)) ,
\end{align}  
and 
\begin{align*}
a_h(\textbf{u}_h,\textbf{v}) :&=  ({\nabla\textbf{u}_h, \nabla\textbf{v}} )-\sum_{E \in \cE^{B}_{h}} \int_{E} \frac{\partial{\textbf{u}_h}}{\partial{\mathbf{n}}} \cdot{\textbf{v}} \ds-\sum_{E \in \cE^{B}_{h}} \int_{E} \frac{\partial{\textbf{v}}}{\partial{\mathbf{n}}} \cdot {\textbf{u}_h  } \ds \nonumber \\&\quad +\sum_{E \in \cE^{B}_{h}} \int_{E} \frac{\zeta}{h_E}{\textbf{u}_h \cdot \textbf{v}} \ds, \nonumber \\ b_h(p_h,\textbf{v}) :&= (p_h,\nabla\cdot\textbf{v})-\sum_{E \in \cE_{h}^{B}}\int_{E}(\textbf{v} \cdot \textbf{n}) p_h \ds, \nonumber \\
S_{h}((\textbf{u}_h,p_h),(\textbf{v},q)):&=S_{si}((\textbf{u}_h,p_h),(\textbf{v},q))+S_{sb}((\textbf{u}_h,p_h),(\textbf{v},q)), \nonumber \\ 
S_{si}((\textbf{u}_h,p_h),(\textbf{v},q)) :&= \sum_{a \in \cV_{h}} \beta_{a} \int_{\cM_{a}}{\kappa_{a}(\nabla \cdot {\textbf{u}_h})}{\kappa_{a}(\nabla \cdot {\textbf{v}})} \dx \nonumber \\&\quad+\sum_{a \in \cV_{h}} \beta_{a} \int_{\cM_{a}}{\kappa_{a}(\nabla{p_h})}{\kappa_{a}(\nabla{q})} \dx, \nonumber \\
S_{sb}((\textbf{u}_h,p_h),(\textbf{v},q)) :&= \sum_{E \in \cE^{B}_{h}} \int_{E} ({\textbf{u}_h \cdot \textbf{n}})({\textbf{v} \cdot \textbf{n}}) \ds,\\
L(\textbf{v}, q) :&= (\textbf{f},\textbf{v}).
\end{align*}
Further, introduce the generalized local projection norm for $\textbf{V}_{h} \times Q_h $  by
\begin{equation}\label{GLP_norm_s}
\tnorm {(\textbf{u}_h,p_h)}^{2} := \Vert \nabla\textbf{u}_h \Vert^{2} + \Vert p_h \Vert^{2} +\sum_{E \in \cE^{B}_{h}} \int_{E} \frac{\zeta}{h_E}{\textbf{u}^{2}_h } \ds+S_{h}( {(\textbf{u}_h,p_h)}, {(\textbf{u}_h,p_h)}).
\end{equation}
\begin{theorem} \label{inf_sup_s_1} Let $\zeta$ be chosen such that $\zeta>\zeta_{0}>0$ with sufficiently large $\zeta_{0}$ and {the discrete bilinear form (\ref{weak_f_cond_1}) satisfies the following inf-sup condition for some positive constant $\nu$, independent of $h$,  }  
\begin{equation*} %\label{inf_sup_s}
 \inf_{(\mathbf{u}_h,p_h) \in \textbf{V}_{h} \times Q_h } \sup_{(\mathbf{v}_h,q_h) \in \textbf{V}_{h} \times Q_h }\frac{B_{h}((\mathbf{u}_h,p_h),(\mathbf{v}_h,q_h))}{ \tnorm{(\mathbf{u}_h,p_h)}  \tnorm{ (\mathbf{v}_h,q_h) }} \geq \nu.
\end{equation*}
\end{theorem}
{\bf Proof.} In order to prove the stability result, it is enough to choose some $(\textbf{v}_h, q_h) \in \textbf{V}_{h} \times Q_h$ for any arbitrary $(\textbf{u}_h, p_h) \in \textbf{V}_{h} \times Q_h,$ such that
\begin{equation*} %\label{inf_sup_s} 
 \sup_{(\textbf{v}_h,q_h) \in \textbf{V}_{h} \times Q_h }\frac{B_{h}((\textbf{u}_h,p_h),(\textbf{v}_h,q_h))}{\tnorm {(\textbf{v}_h,q_h)}} \geq \nu \tnorm{ (\textbf{u}_h,p_h)} > 0.
\end{equation*}
We first consider the bilinear form in (\ref{bilinear_s}) with  $(\textbf{v}_h, q_h) = (\textbf{u}_h, p_h)$.  
\begin{align} \label{s_1}
 B_{h}((\textbf{u}_h,p_h),(\textbf{u}_h,p_h)) &=\norm{\nabla\textbf{u}_h}^{2}  -2\sum_{E \in \cE^{B}_{h}} \int_{E} \frac{\partial{\textbf{u}_h}}{\partial{\mathbf{n}}} \cdot{\textbf{u}_h} \ds +\sum_{E \in \cE^{B}_{h}} \int_{E} \frac{\zeta}{h_E}{\textbf{u}^{2}_h} \ds \nonumber \\ & \qquad + S_{h}((\textbf{u}_h,p_h),(\textbf{u}_h,p_h)).
\end{align}
The second term of (\ref{s_1}) is handled by using the Cauchy-Schwarz inequality and trace inequality (\ref{trace_ineq1})
\begin{align}\label{s_2}
2 \int_{E} \frac{\partial{\textbf{u}_h}}{\partial{\mathbf{n}}} \cdot{\textbf{u}_h} \ds \leq 2\norm{\frac{\partial{\textbf{u}_h}}{\partial{\mathbf{n}}}}_{\rm{L}^{2}(E)} \norm{\textbf{u}_h}_{\rm{L}^{2}(E)} \leq 2h^{-1/2}_{E}\norm{\nabla\textbf{u}_h}_{\rm{L}^{2}(K)} \norm{\textbf{u}_h}_{\rm{L}^{2}(E)}.
\end{align}
The sum of all boundary edges of (\ref{s_2}) and using the Young's inequality we have
\begin{align}\label{s_3}
2\sum_{E \in \cE^{B}_{h}} \int_{E} \frac{\partial{\textbf{u}_h}}{\partial{\mathbf{n}}} \cdot{\textbf{u}_h} \ds \leq \frac{1}{2}\norm{\nabla\textbf{u}_h}^{2}+2C^{2}\sum_{E \in \cE^{B}_{h}} \int_{E} \frac{1}{h_E}{\textbf{u}^{2}_h} \ds. 
\end{align}
Substitution of (\ref{s_3}) to (\ref{s_1}) and the selection of parameter $\zeta>\zeta_{0}=:4C^{2}$ to obtain
\begin{align*}
 B_{h}((\textbf{u}_h,p_h),&(\textbf{u}_h,p_h))\\ &\geq\frac{1}{2}\norm{\nabla\textbf{u}_h}^{2}  +\frac{\zeta-2C^2}{\zeta}\sum_{E \in \cE^{B}_{h}} \int_{E} \frac{\zeta}{h_E}{\textbf{u}^{2}_h} \ds + S_{h}((\textbf{u}_h,p_h),(\textbf{u}_h,p_h)) \\ &\geq\frac{1}{2}\left(\norm{\nabla\textbf{u}_h}^{2}+\sum_{E \in \cE^{B}_{h}} \int_{E} \frac{\zeta}{h_E}{\textbf{u}^{2}_h} \ds + S_{h}((\textbf{u}_h,p_h),(\textbf{u}_h,p_h))\right).
\end{align*}
Note that the selection of parameter $\zeta$ implies that
\begin{align*}
\frac{\zeta-2C^2}{\zeta} \geq \frac{1}{2}.
\end{align*}
 Finally, taking $(\textbf{v}_h, q_h) = (\textbf{z}_h, 0)$ as a test function pair, the bilinear form   
    (\ref{bilinear_s}) becomes
\begin{align}\label{z_h1_1}
B_{h}((\textbf{u}_h,p_h),(\textbf{z}_h, 0)) &= a_{h}(\textbf{u}_h,\textbf{z}_h)-b_{h}(p_h,\textbf{z}_h) + S_{h}((\textbf{u}_h,p_h),(\textbf{z}_h, 0)).
\end{align}
Most of the estimates for the right-hand side terms of (\ref{z_h1_1}) follows from (\ref{z_h1}). Only those estimates that are new or different from (\ref{z_h1}) are discussed here.  Consider the first term of (\ref{z_h1_1})
\begin{align}\label{s_5}
a_{h}(\textbf{u}_h,\textbf{z}_h)&=({\nabla\textbf{u}_h, \nabla\textbf{z}_h} ) -\sum_{E \in \cE^{B}_{h}} \int_{E} \frac{\partial{\textbf{u}_h}}{\partial{\mathbf{n}}} \cdot{\textbf{z}_h} \ds-\sum_{E \in \cE^{B}_{h}} \int_{E} \frac{\partial{\textbf{z}_h}}{\partial{\mathbf{n}}} \cdot {\textbf{u}_h  } \ds \nonumber \\ &\quad +\sum_{E \in \cE^{B}_{h}} \int_{E} \frac{\zeta}{h_E}{\textbf{u}_h \cdot \textbf{z}_h} \ds.
\end{align}
The first term of (\ref{s_5}) is handled by using the Cauchy-Schwarz inequality, (\ref{surjectivity_h1}) and the Young's inequality
\begin{align*}
({\nabla\textbf{u}_h, \nabla\textbf{z}_h} ) \leq \norm{\nabla\textbf{u}_h}\norm{\nabla\textbf{z}_h} \leq C \norm{\nabla\textbf{u}_h}\norm{p_h}  \leq C\norm{\nabla\textbf{u}_h}^{2} +\frac{\norm{p_h}^{2}}{5}.
\end{align*}
Since $\frac{\partial{\textbf{u}_h}}{\partial{\mathbf{n}}}$ is constant on the edge $E$ and $\textbf{z}=0$ on the boundary edges, using the Cauchy-Schwarz inequality, the trace inequality, (\ref{eq_30}) and (\ref{surjectivity_h1})  
\begin{align*}
 \int_{E} \frac{\partial{\textbf{u}_h}}{\partial{\mathbf{n}}} \cdot{\textbf{z}_h} \ds  =\int_{E} \frac{\partial{\textbf{u}_h}}{\partial{\mathbf{n}}} \cdot({\textbf{z}_h}-{\textbf{z}})\ds &\leq C \norm{\frac{\partial{\textbf{u}_h}}{\partial{\mathbf{n}}}}_{\rm{L}^{2}(E)}\norm{{\textbf{z}_h}-{\textbf{z}}}_{\rm{L}^{2}(E)} \\ &\leq C\norm{\nabla\textbf{u}_h}_{\rm{L}^{2}(K)} \norm{\nabla\textbf{z}_h}_{\rm{L}^{2}(K)}.
\end{align*}
and
\begin{align*}
 \int_{E} \frac{\zeta}{h_E}{\textbf{u}_h \cdot \textbf{z}_h} \ds& =\left(\int_{E} \frac{\zeta}{h_E} \textbf{u}^{2}_h  \ds \right)^{\frac{1}{2}}\left(\int_{E} \frac{\zeta}{h_E} ({\textbf{z}_h}-{\textbf{z}})^{2} \ds \right)^{\frac{1}{2}}\\ &\leq C \left(\int_{E} \frac{\zeta}{h_E} \textbf{u}^{2}_h  \ds \right)^{\frac{1}{2}} \norm{\nabla\textbf{z}_h}_{\rm{L}^{2}(K)}.
\end{align*}
The sum of all boundary edges  and using the  (\ref{surjectivity_h1}) and Young's inequality we have
\begin{align*}
\sum_{E \in \cE^{B}_{h}} \int_{E} \frac{\partial{\textbf{u}_h}}{\partial{\mathbf{n}}} \cdot{\textbf{z}_h} \ds &+\sum_{E \in \cE^{B}_{h}} \int_{E} \frac{\zeta}{h_E}{\textbf{u}_h \cdot \textbf{z}_h} \ds \\ &\leq C\left(\norm{\nabla\textbf{u}_h}^{2}+ \sum_{E \in \cE^{B}_{h}}\int_{E} \frac{\zeta}{h_E} \textbf{u}^{2}_h  \ds  \right) +\frac{\norm{p_h}^{2}}{5}.
\end{align*}
The second term of (\ref{s_5}) is handled as
\begin{align*}
\sum_{E \in \cE^{B}_{h}} \int_{E} \frac{\partial{\textbf{z}_h}}{\partial{\mathbf{n}}} \cdot{\textbf{u}_h} \ds &\leq \norm{\nabla\textbf{z}_h}^{2}+C\sum_{E \in \cE^{B}_{h}} \int_{E} \frac{\zeta}{h_E} {\textbf{u}^{2}_h} \ds  \leq \frac{\norm{p_h}^{2}}{5}+C\sum_{E \in \cE^{B}_{h}} \int_{E} \frac{\zeta}{h_E} {\textbf{u}^{2}_h} \ds.
\end{align*}
Put together (\ref{z_h1_1}) leads to
\begin{align*}
B_{h}((\textbf{u}_h,&p_h),(\textbf{z}_h, 0)) \\ &\geq \frac{1}{2}\norm{p_h}^{2}-C\left(\norm{\nabla\textbf{u}_h}^{2}+ \sum_{E \in \cE^{B}_{h}}\int_{E} \frac{\zeta}{h_E} \textbf{u}^{2}_h  \ds +S_{h}((\textbf{u}_h,p_h),(u_h,p_h)) \right).
\end{align*}
The final selection of $(\textbf{v}_h,q_h)$ is
\begin{align*}
(\textbf{v}_h, q_h) = (\textbf{u}_h, p_h) +\frac{1}{{C}+1}(\textbf{z}_h, 0),
\end{align*}
here $I_h$ is  defined in ({\ref{eq_30}}). %(\ref{eq_25a}) and (\ref{eq_24a}) respectively.
Finally, 
%(\ref{eq_26ab}), (\ref{eq_25a}) and (\ref{eq_24a}) and 
 rest of the proof follows in a similar way as in Theorem \ref{inf_sup_d_1}.
\begin{theorem} \label{th_012}
 Let $(\mathbf{u}, p) \in [{\rm{H}^{2}}(\Omega)]^{2}\times  \rm{L}^{2}_{0} \bigcap \rm{H}^{1}(\Omega)\  and \ \ (\mathbf{u}_h, {p}_h) \in \textbf{V}_h \times\  Q_h$ be the solutions to (\ref{weak_f_d}) and (\ref{weak_f_cond_1}), respectively. Suppose $\beta_a = \beta{h_a};$ for some $\beta > 0.$ Then it holds  
  \begin{equation*}
  \tnorm {(\mathbf{u}-\mathbf{u}_h, p-p_h)} \leq C \left(\norm{h_{\T} \mathbf{u}}_{2} + \norm{ h_{\T} p}_{1} \right).
  \end{equation*}
  \end{theorem}
{ \bf Proof.} Identical to the proof of Theorem  \ref{th_01}.
%%%%%%%%%%%%%%%%%%%%%%%%%%%%%%%%%%%%%%%
\section {Numerical Results} \label{computation}
\begin{figure}[!ht]
	\centerline{%
		\begin{tabular}{cc}
			\hspace{0cm}
			\resizebox*{7cm}{!}{\includegraphics{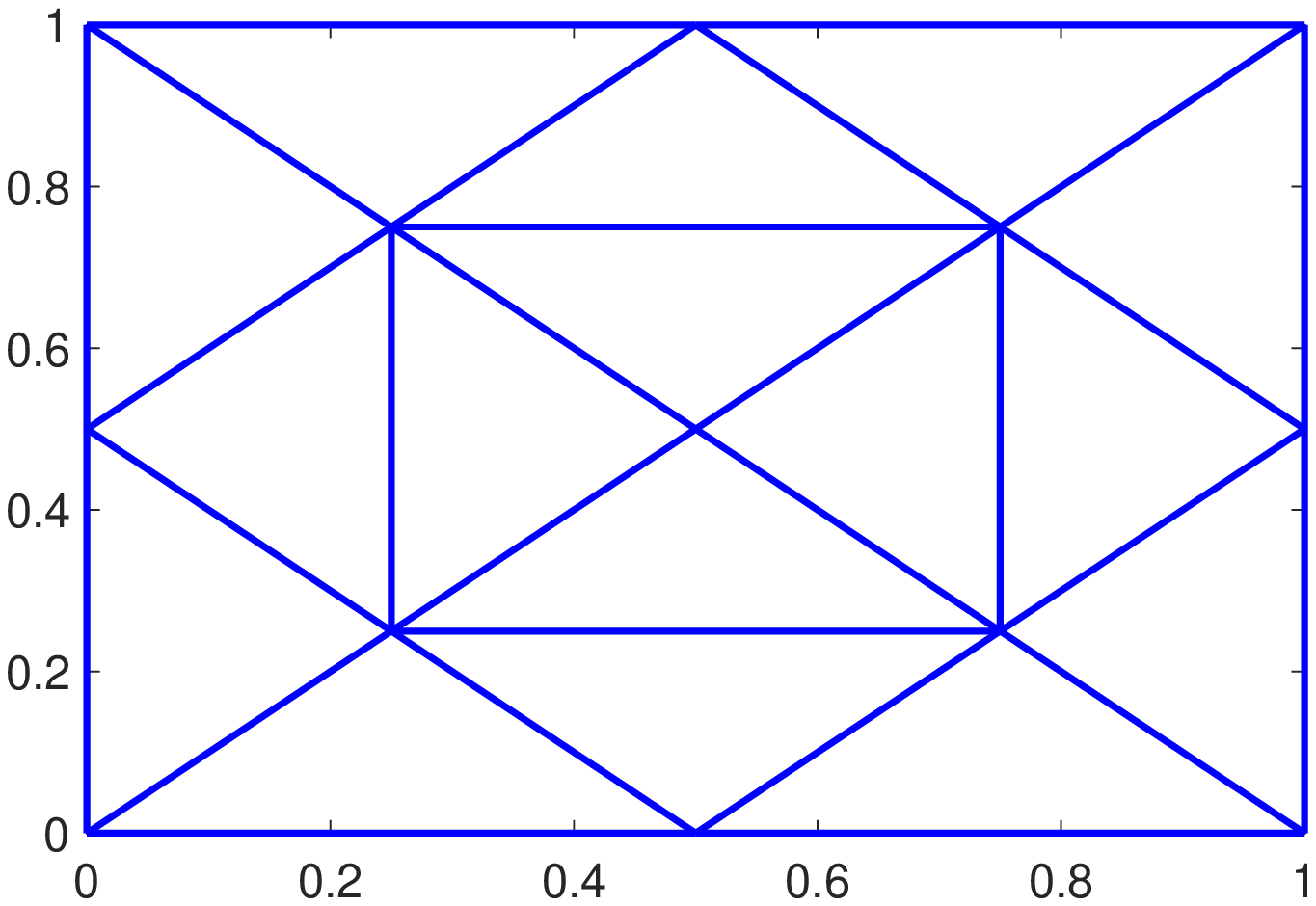}}%
			&\hspace{-1cm}
			\resizebox*{7cm}{!}{\includegraphics{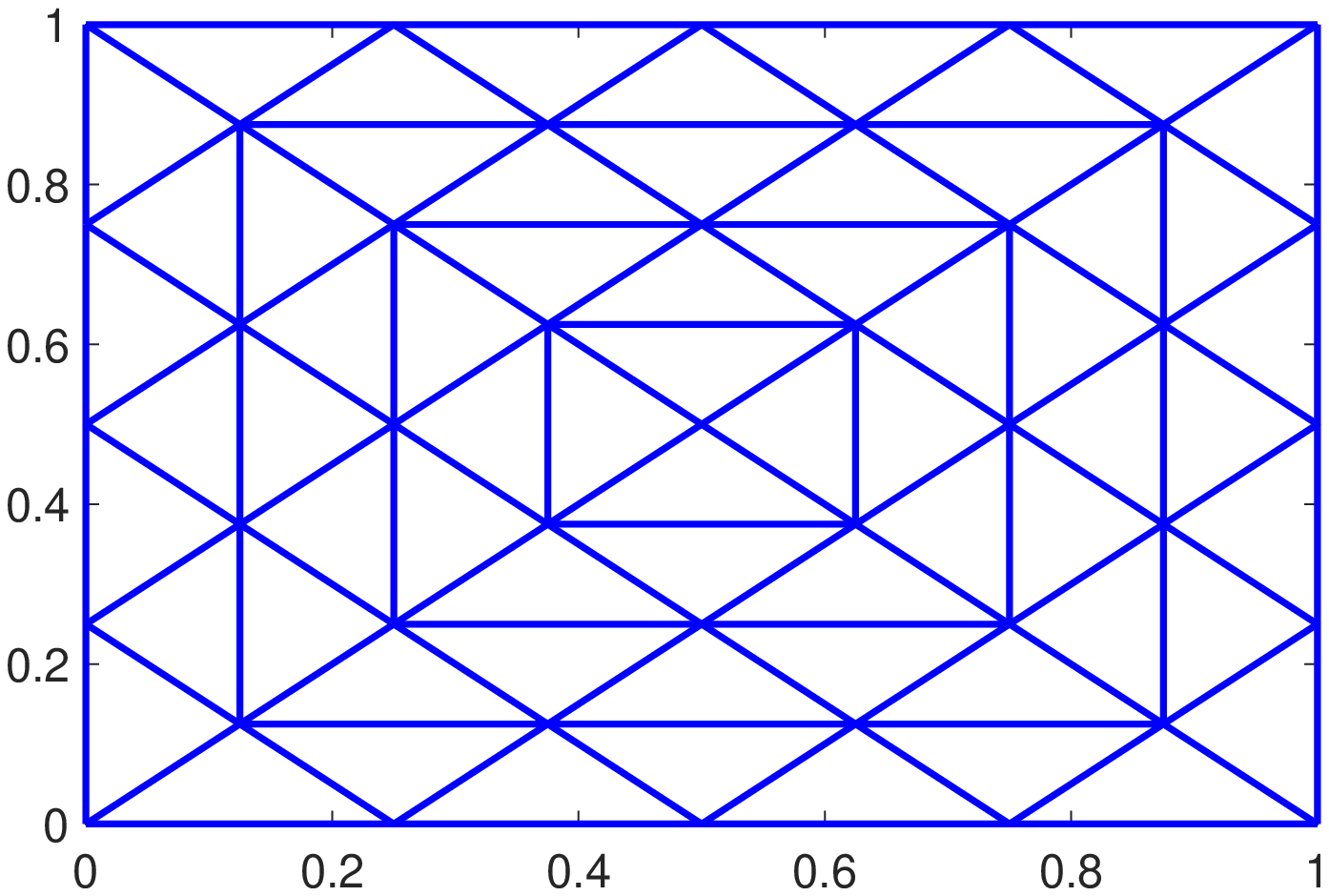}}%
			\\
			%{(a)\ Initial Mesh} \hspace{-1cm}&{(b) \ Refinement of mesh}
		\end{tabular}
	} \caption{{\label{mesh}{Triangulations used for computations in section \ref{computation}}}}
\end{figure}

In this section, we present an array of numerical results to support the derived theoretical estimates and to illustrate the robustness of the proposed scheme.
 Numerical solutions of all test examples are computed on an hierarchy  of uniformly refined triangular meshes having 16, 64, 256, 1024, and 4096 cells, respectively, see   Figure~\ref{mesh} for the initial and an uniformly refined mesh {of 16 triangles and 64 triangles, respectively.}   
 
 \subsection*{A. Darcy flow problem} \label{ex1}
\begin{figure}[ht!]
\centerline{\includegraphics[width=6.0cm]{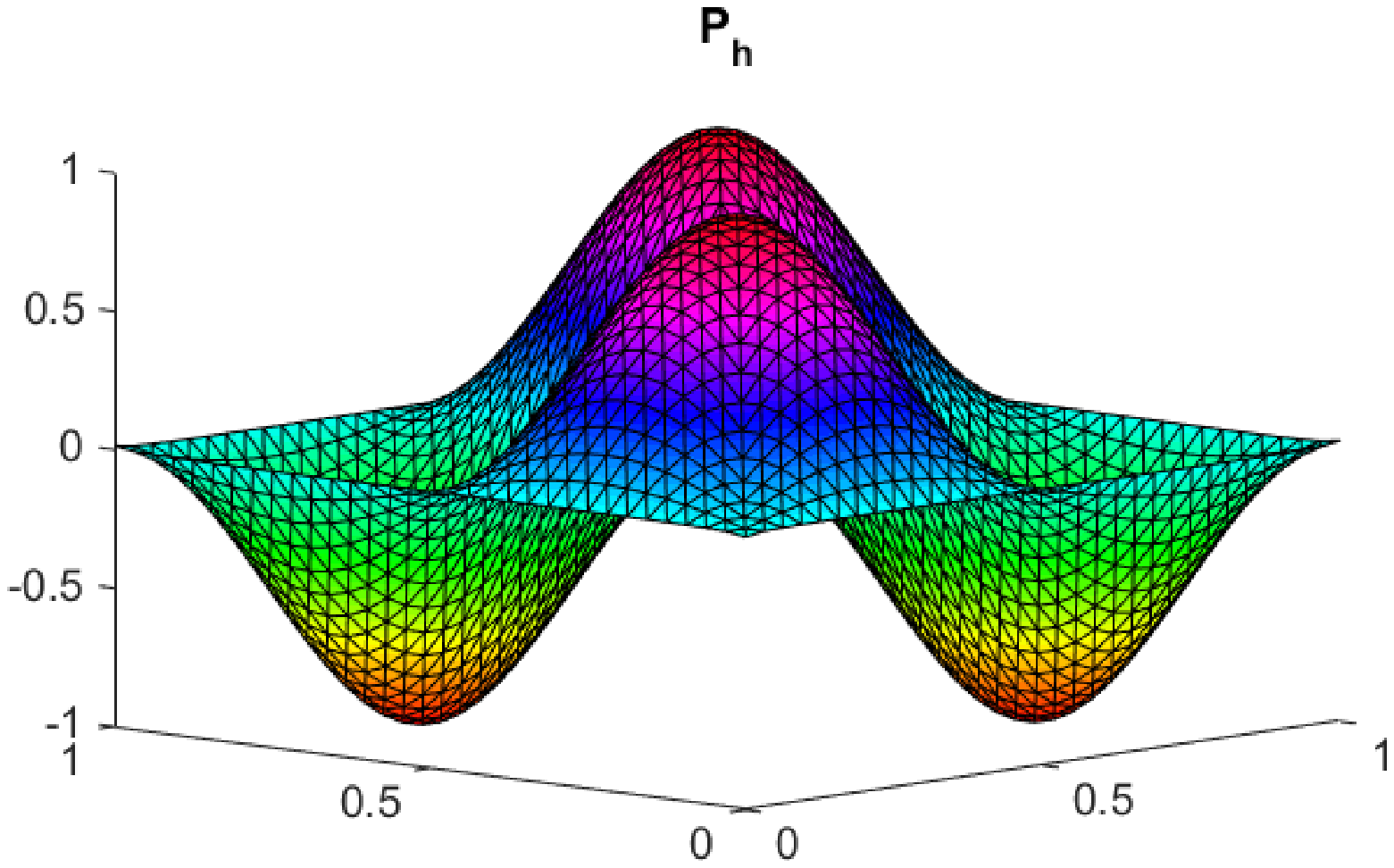}\hspace*{1cm}\includegraphics[width=6.0cm]{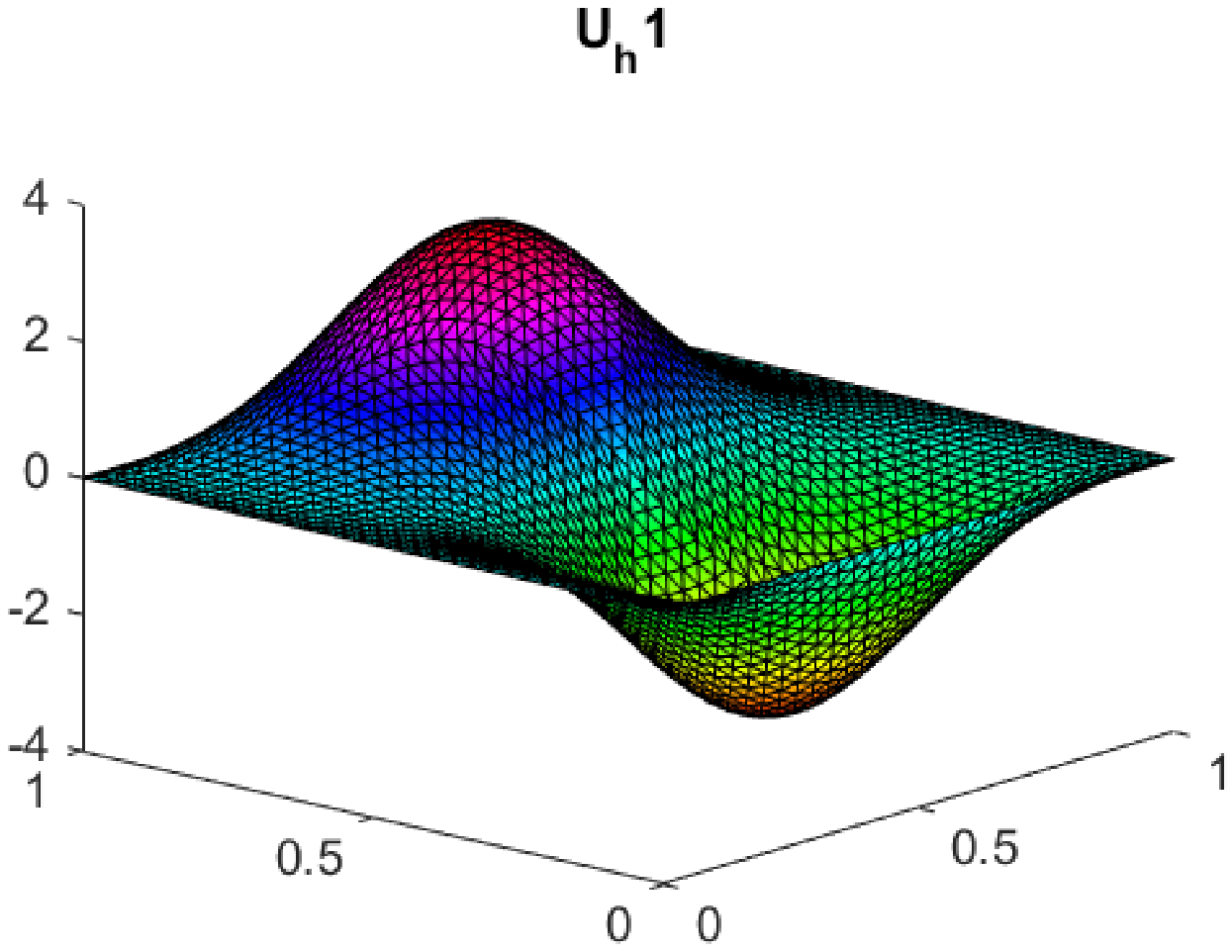}}
\centerline{\includegraphics[width=6.0cm]{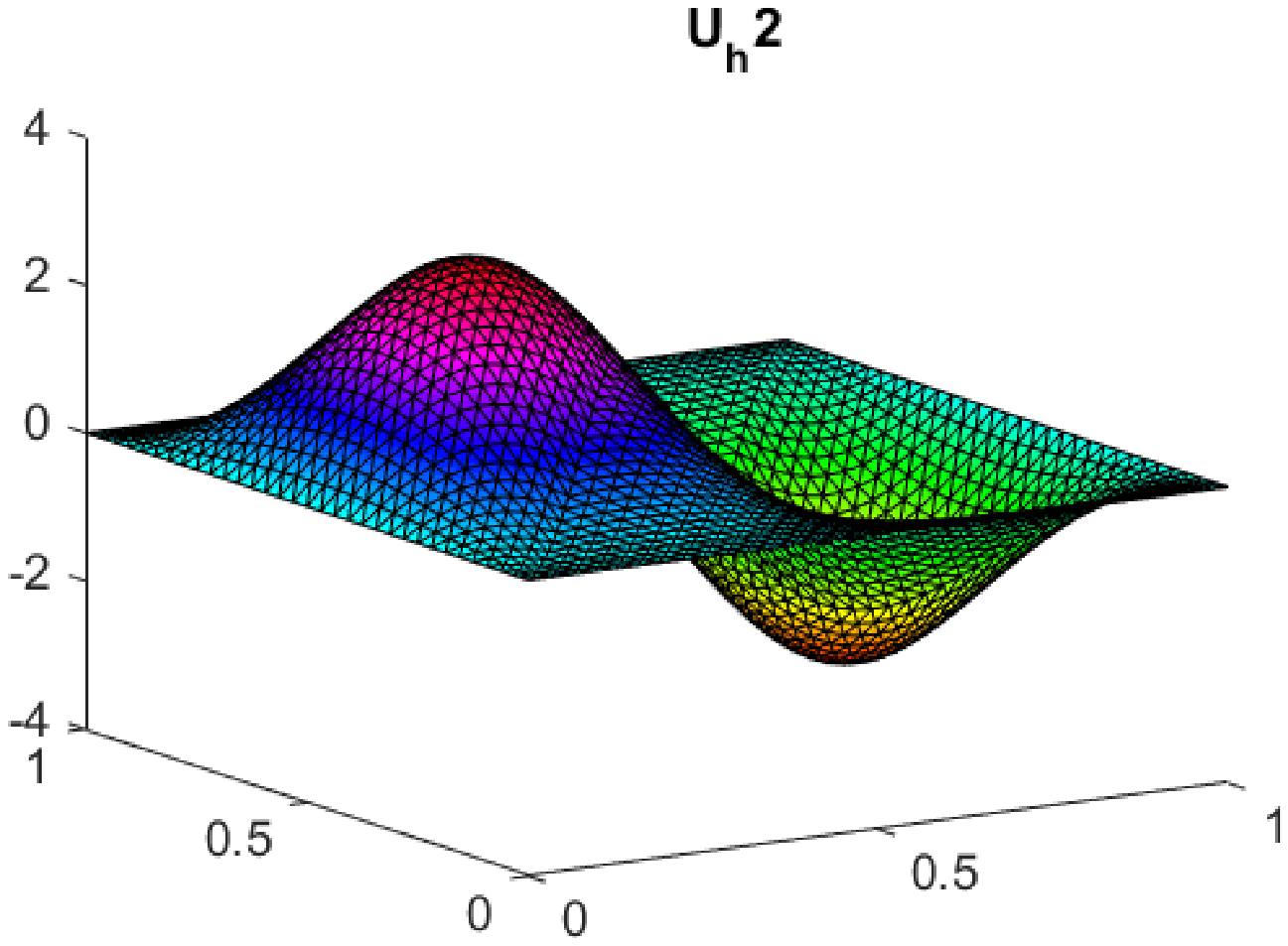}\hspace*{0cm}\includegraphics[width=7.5cm]{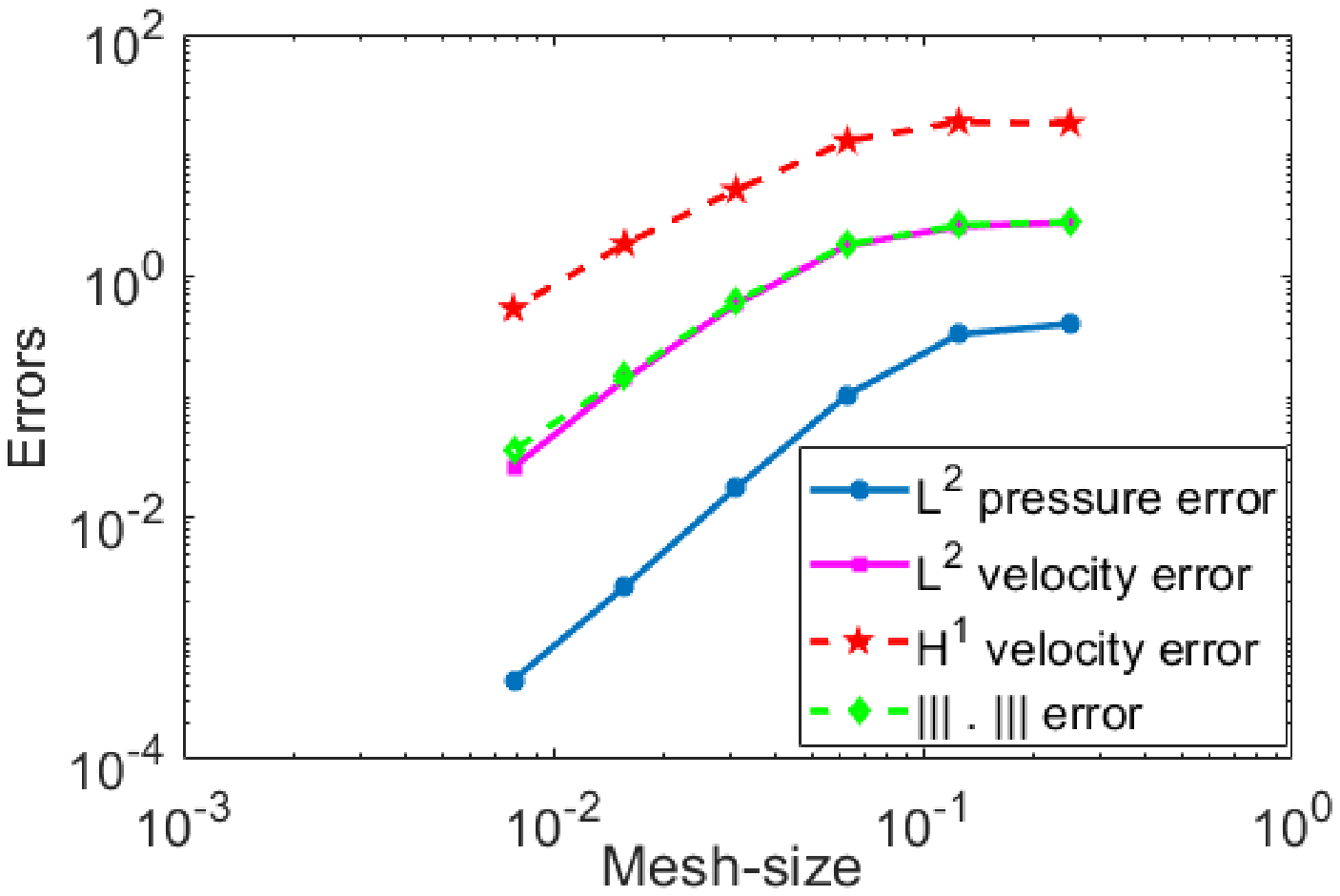}}
\caption{{ GLPS discrete solution $(\textbf{u}_h , p_h)$, and convergence plot of Darcy problem.} }\label{darcy_fig}
\end{figure}

 We consider the model problem \eqref{problem1}  in $\Omega=(0, 1)^2$  with  a given exact solution 
 \[ 
\textbf{u}(x,y)=(-\pi \sin(2 \pi y) \sin^{2}(\pi x),\pi \sin(2\pi x)  \quad \text{and} \quad  p(x,y)=\sin(2\pi x)\sin(2\pi y)
\]
and set a stabilization parameter $\beta_a = \beta h_a$, $\beta=10$. The solution is approximated with the equal-order interpolation spaces $\mathbf{P}_{1}^{c}/\mathbf{P}_{1}^{c}$ using GLPS finite element formulation (\ref{weak_f_cond}). Although the velocity and pressure approximation spaces are not {\textit{inf-sup}} stable for the Darcy problem, the GLP stabilization arrests the oscillations  effectively.  
Figure \ref{darcy_fig} shows the $\mathbf{P}_{1}^{c}/\mathbf{P}_{1}^{c}$ approximations with GLP stabilized finite element solutions   {at the mesh-size} 0.0078.    {The errors are computed in   $ \rm{L}^{2} $- norm,  $\rm{H}^{1}$-seminorm  and $\tnorm{\cdot}$ stabilized  norm}. The computed errors with the $ \rm{L}^2-$norm and  $\rm{H}^{1}-$seminorm are presented in {Table~\ref{darcy_table1}, whereas Table~\ref{triple_darcy_stoke} presents the errors measured in GLP stabilized norm as defined in \eqref{GLP_norm_d}.  
} We can observe a second-order convergence in $ \rm{L}^{2}$-norm, a first-order convergence in $\rm{H}^{1}$-seminorm  {and $\mathcal{O}(h^{3/2})$ convergence in $\tnorm{\cdot}$}. Also, the last plot of Figure \ref{darcy_fig} shows the convergence behavior of $\mathbf{P}_{1}^{c}/\mathbf{P}^{c}_{1}$ approximation of Darcy equations with respect to $\rm{L}^{2}$-norm, $\rm{H}^{1}$-seminorm and  {the GLP stabilized norm}. These numerical results support the estimates derived in the previous section.

%%%%%%%%%%%%%%%%%%%%%%%%%%%%%%%%%%%%%%%%%%%%

%%%%%%%%%%%%%%%%%%%%%%%%%%
 \begin{table}[h!]
 	\begin{center}
 		\caption{{{Darcy problem: Errors and convergence orders. }} }\label{darcy_table1}\begin{tabular}{|l|c|c|c|c|c|r|} % <-- Changed to S here.
 			%\textbf{Value 1} & \textbf{Value 2} & \textbf{Value 3} & \textbf{Value 3}& \textbf{Value 3}\\
 			\hline
 			$\mbox{Mesh-size}$ & $ \norm{\mathbf{u}-\mathbf{u}_h}$ & $\mbox{Order}$& $ |{\nabla(\mathbf{u}-\mathbf{u}_h)}|$ & $\mbox{Order}$ & $ \norm{p-p_h}$& $\mbox{Order}$\\ [.75ex]
 			\hline  
%  			1/4 &     2.7610 &  -  &     18.4120      &        -          &      0.4006 &  - \\ [.75ex]
%  			1/8    &   2.5773  &  0.0994   &   18.8393   &      -0.0331    &       0.3289   &   0.2848\\ [.75ex]
 			1/16    & 1.7949 &  -   &    13.1479  &          -    &        0.1040   &   - \\ [.75ex]
 			1/32  &   0.5847   &  1.6182 &    5.1579  &         1.3500    &      0.0177 &    2.5549\\ [.75ex]
 			1/64   &   0.1395 &   2.0669 &   1.8466  &          1.4819      &     0.0027 &   2.7185\\ [.75ex]
 			1/128 &  0.0262 &   2.4128 &   0.5405 &         1.7724      &      0.0005 &   2.5674
\\ [.75ex]
 			\hline
 		\end{tabular}
 	\end{center}
 \end{table}

%%%%%%%%%%%%%%%%%%%%%%%%%%%%%

\subsection*{B. Stokes flow problem} \label{ex2}
\begin{figure}[ht!]
\centerline{\includegraphics[width=6.0cm]{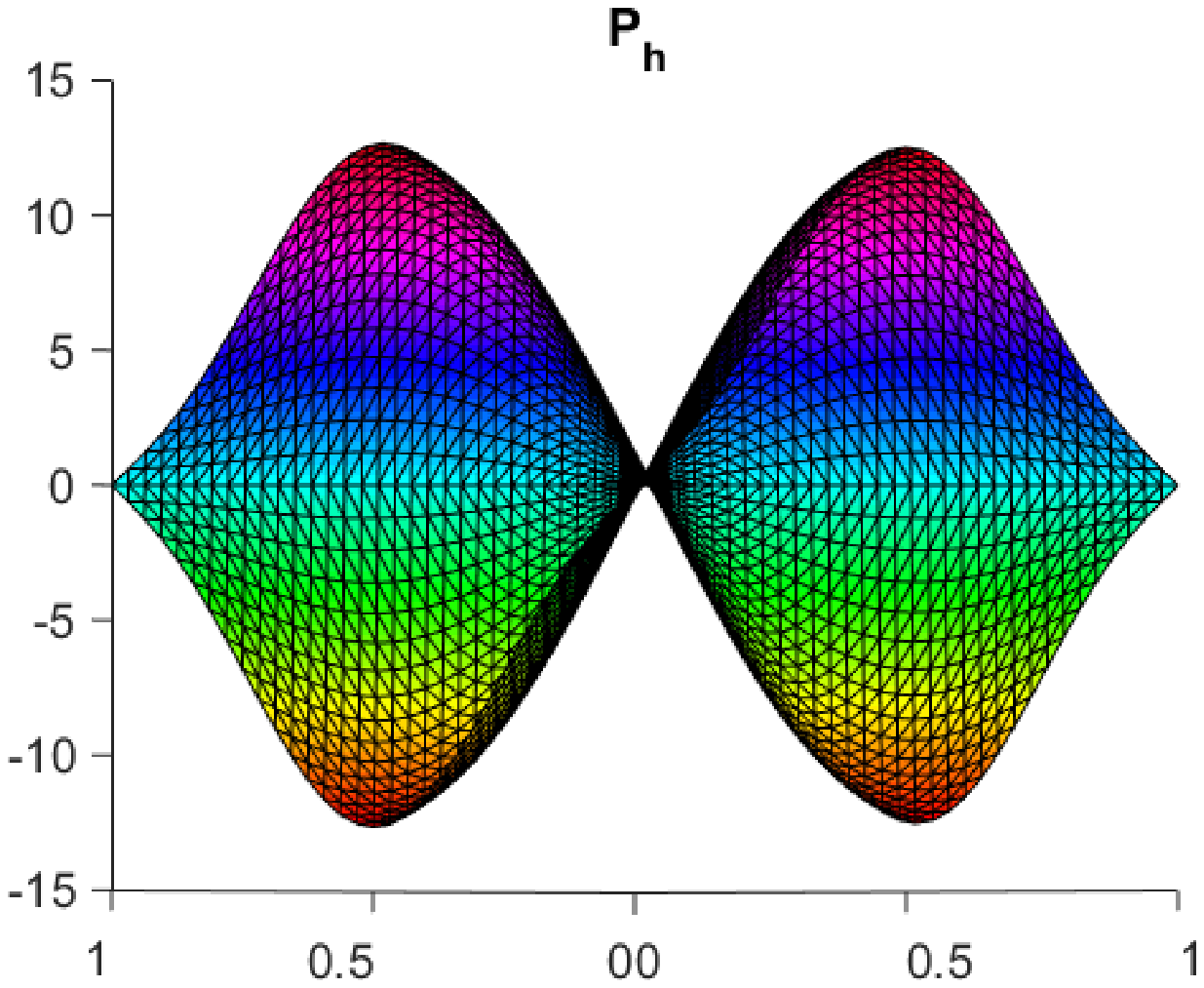}\hspace*{1cm}\includegraphics[width=6.0cm]{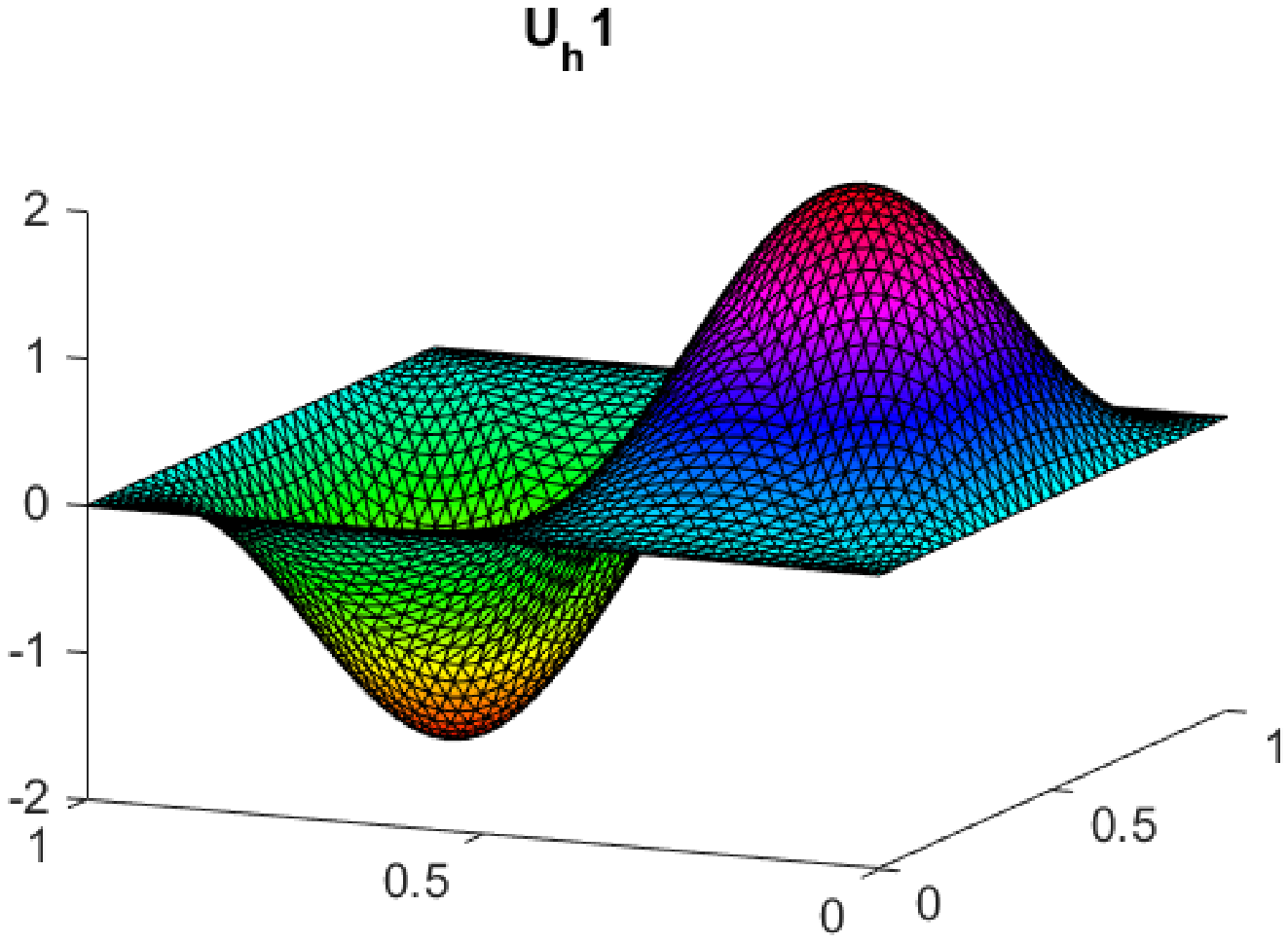}}
\centerline{\includegraphics[width=6.0cm]{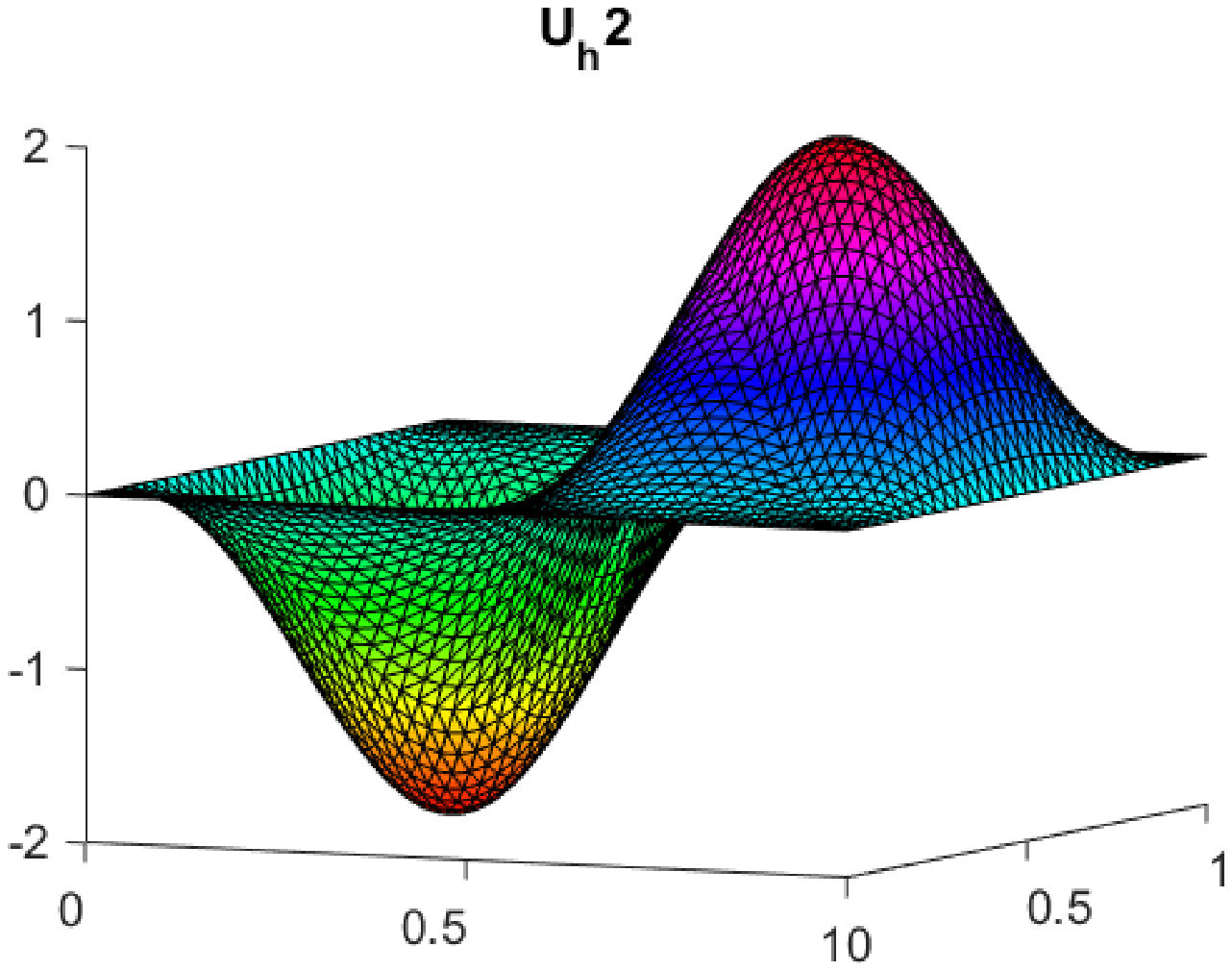}\hspace*{0cm}\includegraphics[width=7.5cm]{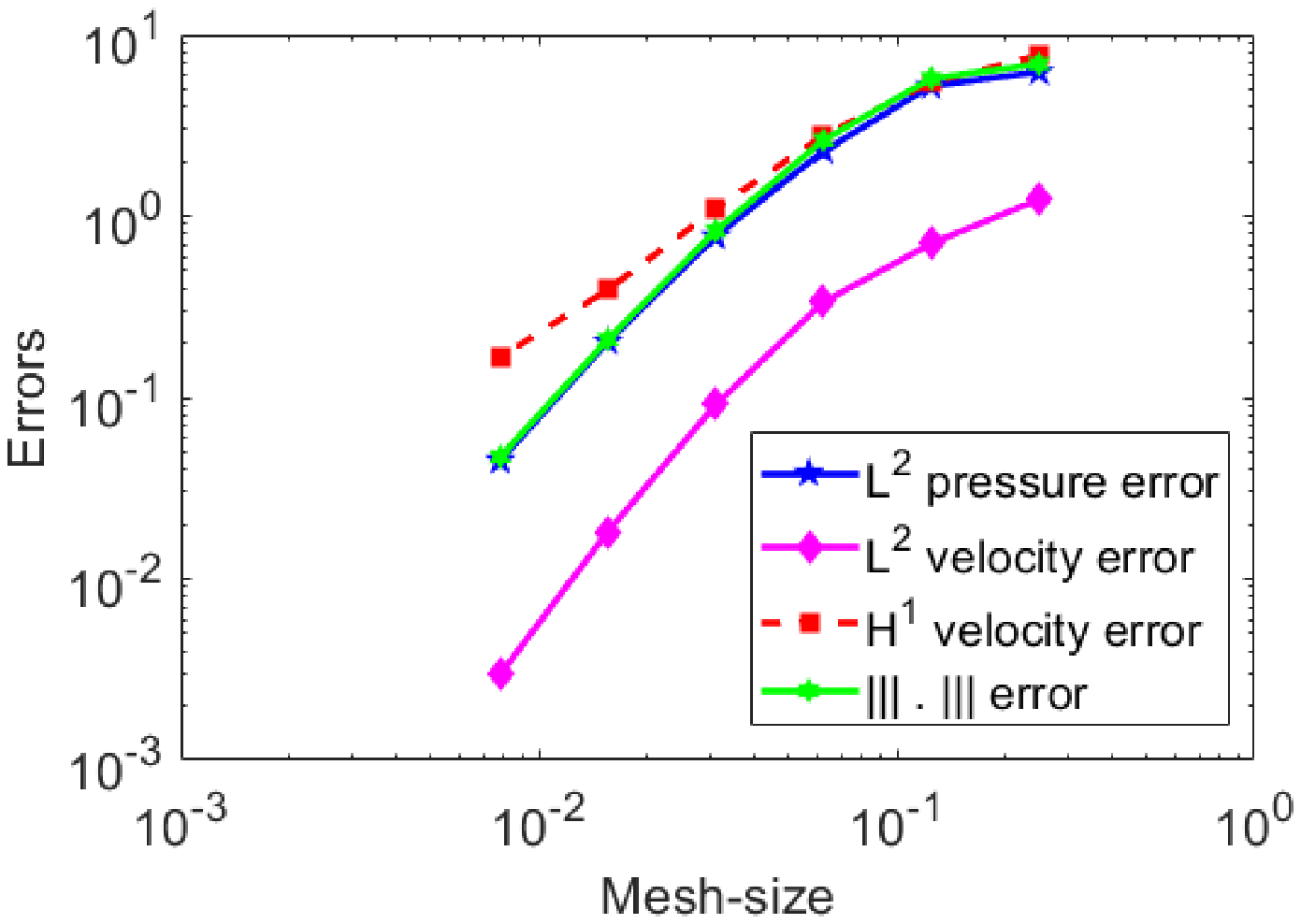}}
\caption{{GLPS discrete solution $(\textbf{u}_h , p_h)$  and convergence plot of Stokes problem.} }\label{stoke_fig}
\end{figure}

In order to demonstrate the robustness of the method, we consider the Stokes problem as the second numerical test example.  We consider the model problem \eqref{stoke}  in $\Omega=(0, 1)^2$  with  a given exact solution ${ \textbf{\rm \textbf{u}}}(x,y)=(-\cos(2 \pi x) \sin(2 \pi y)+\sin(2 \pi y), \sin(2 \pi x) \cos(2 \pi y)-\sin(2 \pi x))$ and $ p(x,y)=2 \pi (\cos(2 \pi y)-\cos(2 \pi x)).$
The stabilization parameters for the discrete variational formulation (\ref{bilinear_s}) is chosen as $\beta_a = h_a$ with $\beta=1$ and $\zeta=2$.  The equal-order interpolation spaces, $\mathbf{P}_{1}^{c}/\mathbf{P}_{1}^{c}$, are used to approximate the velocity and pressure approximation. The generalized local projection stabilized finite element scheme overcomes the space incompatibility issue and improves the pressure's approximation. Figure \ref{stoke_fig} displays the $\mathbf{P}_{1}^{c}/\mathbf{P}_{1}^{c}$ stabilized solution    {at the mesh-size} 0.0078. The quantitative and qualitative errors and the order of convergence obtained with $\mathbf{P}_{1}^{c}/\mathbf{P}_{1}^{c}$ finite element approximations are summarized in Table \ref{stoke_table1_3},  {Table \ref{triple_darcy_stoke} } and in  the last plot of Figure~\ref{stoke_fig}. 
Desired convergence rates, $i.e.$, second-order $\rm{L}^{2}$-errors in velocity and pressure, and first-order   $\rm{H}^{1}-$approximation error in velocity, are demonstrated.

%%%%%%%%%%%%%%%%%%%%%%%%%%%%%%%%%%%%%%%%%%%%%

 \begin{table}[h!]
 	\begin{center}
 		\caption{{{Stokes problem: Errors and convergence orders. }} }
 		\label{stoke_table1_3}
 		\begin{tabular}{|l|c|c|c|c|c|r|} % <-- Changed to S here.
 			%\textbf{Value 1} & \textbf{Value 2} & \textbf{Value 3} & \textbf{Value 3}& \textbf{Value 3}\\
 			\hline
 			$\mbox{Mesh-size}$ & $ \norm{\mathbf{u}-\mathbf{u}_h}$ & $\mbox{Order}$& $ |{\nabla(\mathbf{u}-\mathbf{u}_h)}|$ & $\mbox{Order}$ & $ \norm{p-p_h}$& $\mbox{Order}$\\ [.75ex]
 			\hline  
%  			1/4 &    1.2345 &  -  &     7.6325     &        -          &      6.2216 &  - \\ [.75ex]
%  			1/8    &    0.7050 &   0.8083   &    5.4559     &      0.4843     &      5.2113   &   0.2556 \\ [.75ex]
 			1/16    & 0.3360  & -  &    2.7865   &          -    &        2.2824  &    -\\ [.75ex]
 			1/32  & 0.0912  & 1.8806   &    1.0935  &          1.3495     &       0.7719  &    1.5641\\ [.75ex]
 			1/64   &   0.0182  &   2.3237 &   0.3933  &          1.4754      &     0.2027  &   1.9289\\ [.75ex]
 			1/128 &  0.0030 &    2.6161 &    0.1667 &        1.2384&          0.0448
 &   2.1764
\\ [.75ex]
 			\hline
 		\end{tabular}
 	\end{center}
 \end{table}

\begin{table}[h!]
  \begin{center}
  \caption{{Error and convergence orders with respect to  $\tnorm{\cdot}$}} \label{triple_darcy_stoke}
    \begin{tabular}{|l|c|c|c|c|c|c|r|} % <-- Changed to S here.
      %\textbf{Value 1} & \textbf{Value 2} & \textbf{Value 3} & \textbf{Value 3}& \textbf{Value 3}\\
           \hline  
     &     $\mbox{Mesh-size} \ h$ &  1/4  &    1/8     &       1/16         &       1/32 & 1/64&1/128 \\ [1ex]
                           Darcy flow                     &   $\tnorm{\cdot}$ & 2.7966   &    2.6309     &     1.8345  &      0.6146  &   0.1470& 0.0363\\ [1ex]
                                     & Order & -  &    0.0881   &         0.5202    &        1.5777 &  2.0639& 2.0164\\ [1ex]
    \hline 
    &  $\mbox{Mesh-size} \ h$  & 1/4  &    1/8  &          1/16    &      1/32 &   1/64& 1/128\\ [1ex]
                      Stokes flow        &   $\tnorm{\cdot}$ &   6.8727 &   5.7263  &         2.6194      &      0.8297  &   0.2087  &0.0470\\ [1ex]
                            &  Order &   - &    0.2633 &         1.1284        &     1.6585 &   1.9910 &2.1508\\ [1ex]
     \hline
    \end{tabular}
  \end{center}
   \end{table}

%%%%%%%%%%%%%%%%%%%%%%%%%%%%
%%%%%%%%%%%%%%%%%

%%%%%%%%%%%%%%%%%%%%%%%%%%%%%%%%%%%%%%%%
%%%%%%%%%%%%%%%%%%%%%%%%%%%%%%%%%%
%%%%%%%%%%%%%%%%%%%%%%%%%%%%%%%%
%%%%%%%%%%%%%%%%%%%%%%%%%%%%%%%%%%%%%
\noindent
\bigskip

%%%%%%%%%%%%%%%%%%%%%%%%%%%%%%%%%

%%%%%%%%%%%%%%%%%%%%%%%%%%%%%%%%%%%
\section{ Conclusions} 
In this article, a generalized local projection stabilized (GLPS) conforming finite element scheme  for Darcy flow and Stokes problems with equal-order interpolation spaces $(\mathbf{P}_{1}^{c}/\mathbf{P}^{c}_{1})$, is proposed and analyzed.  GLPS   allows to use projection spaces on overlapping sets and avoids the need of a two-level mesh or an enrichment of finite element space. 
 The partition of unity of the basis functions together with $\rm{L}^{2}$-orthogonal projection properties is used in deriving the stability and convergence estimates.  Further, a robust  {\it a priori} error analysis is derived for both problems. An array of numerical experiments are presented to support the derived estimates and to demonstrate the efficiency of the proposed scheme in suppressing oscillations without compromising the order of convergence.

\section*{Acknowledgments}  {The work of first author was supported in part by the department of National Mathematics Initiative at IISc Bangalore and the Tata Trusts Travel Grant (ODAA/INT/19/189).}  The second author acknowledges the partially support of Science and Engineering Research Board (SERB) with the grant EMR/2016/003412.

\section*{References}

\bibliographystyle{plain}
\bibliography{bibtexexample}
 
% \bibliographystyle{plain}
% \bibliography{bibtexexample}

\end{document}